\documentclass[12pt]{article}
\usepackage[german]{babel}
\usepackage{amssymb}
\usepackage{amsmath}
\usepackage{graphicx}
\def\noi{\noindent}

\def\IN{\mathbb N}
\def\IZ{\mathbb Z}
\def\IR{\mathbb R}
\def\IC{\mathbb C}
\def\IQ{\mathbb Q}

\def\IA{\mathbb A}
\def\II{\mathbb I}
\def\ik{\Bbbk}

\def\an{\mathrm{an}}
\def\exp{\mathrm{exp}}

\def\vs{\vspace}
\def\ma{\mathcal}
\def\Re{\mathrm{Re}}
\def\Im{\mathrm{Im}}

\begin{document}
	
	\begin{center}
		{\bf \Large Periods, Power Series, and Integrated Algebraic Numbers}
	\end{center}

	\vspace{0.5cm} \centerline{Tobias Kaiser}
	
	\vspace{0.7cm}
	\begin{center}
		\begin{minipage}[t]{10cm}\scriptsize{{\bf Abstract.}
				Periods are defined as integrals of semialgebraic functions defined over the rationals. Periods form a countable ring not much is known about. Examples are given by taking the antiderivative of a power series which is algebraic over the polynomial ring over the rationals and evaluate it at a rational number. We follow this path and close these algebraic power series under taking iterated antiderivatives and nearby algebraic and geometric operations. We obtain a system of rings of power series whose coefficients form a countable real closed field. Using techniques from o-minimality we are able to show that every period belongs to this field. In the setting of o-minimality we define exponential integrated algebraic numbers and show that exponential periods and the Euler constant is an exponential integrated algebraic number. Hence they are a good candiate for a natural number system extending the period ring and containing important mathematical constants.}
		\end{minipage}
	\end{center}
	
	\normalsize
	
	\vs{1cm}
	{\bf 1. Introduction and Main Results}
	
	\vs{0.5cm}
	Algebraic {\bf Periods} are given by an absolutely convergent integral of rational functions, over domains in euclidean spaces given by polynomial inequalities with rational coefficients. Equivalently they are given by an absolutely convergent integral over functions on the whole euclidean space that are semialgebraic over the rationals.
	Kontsevich and Zagier [25] provide a very nice introduction to periods (see also M\"uller-Stach [29]). Note that we consider here real periods defined via sets and functions in euclidean space. This is no restriction since one can separate real and imaginary part (using the canonical identification of $\IC$ with $\IR^2$).
	Periods are a fundamental object in number theory and arithmetic geometry. They arise in bridging algebraic de Rham cohomology and singular cohomology and much more. 
	We refer to Huber and M\"uller-Stach [15] and Huber and W\"ustholz [16] for periods and arithmetic geometry.
	Periods form a countable ring extending the field of algebraic numbers by transcendental numbers. For example, $\pi$ is a period as well as the logarithm of algebraic numbers and values of the Riemann zeta function at natural numbers. It is not known whether the Euler number $e$ or the Euler constant $\gamma$ are periods. So far there are no natural examples of non-periods known (see the end of [15]).
	
	Periods are defined by interaction of algebraic geometry and analysis via integration. Since integration is difficult, periods are so far not very well understood.

	\rule{14cm}{0.01cm}
	
	
	The first goal of this paper is to provide a systematic construction of a countable number system containing the periods where only formal integration is needed.
	
	This approach starts with the following examples for obtaining periods.
	Given a real power series which is algebraic over the rational polynomial ring, take its formal antiderivative and evaluate this power series at a rational point contained in its domain of convergence.
	
	\vs{0.5cm}
	{\bf Examples}
	
	\begin{itemize}
		\item[(1)] Let 
		$$L(X)=\sum_{p=1}^\infty \frac{(-1)^{p+1}}{p} X^p$$ 	
		be the {\bf logarithmic series}. Its radius of convergence is one. 
		We have
		$$\log\Big(\frac{1}{2}\Big)=\int_1^{1/2}\frac{dt}{t}=L\Big(-\frac{1}{2}\Big).$$
		\item[(2)] Let 
		$$A(X)=\sum_{p=0}^\infty \frac{(2p)!}{(2^pp!)^2}\frac{X^{2p+1}}{2p+1}$$ 
		be the {\bf arcsine series}. Its radius of convergence is one.
		We have 
		$$\frac{\pi}{6}=\arcsin\Big(\frac{1}{2}\Big)
		=\int_0^{1/2}\frac{dt}{\sqrt{1-t^2}}=A\Big(\frac{1}{2}\Big).$$
	\end{itemize}
	
	\vs{0.2cm}
	We follow this approach. We start with the power series which are algebraic over the rational polynomial ring. We call them simply {\bf algebraic power series}. We close them under taking formal antiderivatives and nearby algebraic and geometric operations.
	
	\vs{0.5cm}
	{\bf Main Definition}
	
	\vs{0.1cm}
	We consider the smallest class of rings of power series which contain the algebraic series and are closed under taking 
	\begin{itemize}
		\item[(a)] {\bf Reciprocals},
		\item[(b)] {\bf Rational Polynomial Substitutions},
		\item[(c)] {\bf Rational Translations},
		\item[(d)] {\bf Formal Antiderivatives}.
	\end{itemize}
	
	\newpage
	We explain these terms:
	
	\begin{itemize}
		\item[(a)] {\bf Reciprocal:} 
		The reciprocal is defined for power series which do not vanish at the origin. For example, the reciprocal of $1-X$ is the {\bf geometric series} $G(X)=\sum_{p=0}^\infty X^p$.
		\item[(b)] {\bf Rational Polynomial Substitution:}
		We replace the coordinates by rational polynomials vanishing at the origin.
		For example $G(X_1+X_2)=\sum_{p,q=0}^\infty\binom{p+q}{p}X_1^pX_2^q$.
		\item[(c)] {\bf Rational Translation:} 
		Given a convergent power series and a point in its domain of convergence with rational components we develop the series at this point. For example $G(X+1/2)=\sum_{p=0}^\infty 2^{p+1}X^p$.
		\item[(d)] {\bf Formal Antiderivative:} 
		Given a power series we take the formal antiderivative with respect to an arbitrary variable. For example the formal antiderivative of $G(X)$ is given by $\sum_{p=0}^\infty X^{p+1}/(p+1)$.
	\end{itemize}
	
	\vs{0.2cm}
	We call the power series obtained by the main definition the {\bf integrated algebraic power series.} (Note that we have used in [19] the same name in a different set up. There algebraic means algebraicity over the real polynomial ring, here algebraicity over the rational polynomial ring.) 
	
	\vs{0.5cm}
	We obtain the following results. 
	
	\vs{0.5cm}
	{\bf Theorem A}
	
	\vs{0.1cm}
	{\it The coefficients of the integrated algebraic power series form a countable real closed field.}
	
	\vs{0.5cm}
	We call this the {\bf field of integrated algebraic numbers}.
	
	\vs{0.5cm}
	The periods are contained in the established number system.
	
	\vs{0.5cm}
	{\bf Theorem B}
	
	\vs{0.1cm}
	{\it A period is an integrated algebraic number.}
	
	\vs{0.5cm}
	Note that Tent and Ziegler [32] also have constructed a (presumably much larger) countable real closed field containing the periods. There the construction is from a completely different point of view, namely from computability, whereas we work in a tame geometric-analytic setting.
	
	Let us discuss the results here. Periods are defined as integrals of functions which are semialgebraic over the field of rationals. We capture periods as coefficients of integrated algebraic series. Their definition starts with the power series algebraic series which are algebraic over the rational polynomial ring. Hence we start with an algebraic concept. The operations of taking reciprocals, rational substitutions and formal antidervatives are purely formal. The analysis is completely contained in evaluating power series at rational points. 
	
	\vs{0.5cm}
	To establish Theorem A we work with systems of convergent power series with various properties. A central one is given by Weierstra\ss $ $ division and preparation. Denef and Lipshitz [7] have introduced formal Weierstra\ss $ $ systems on fields of characteristic $0$. Van den Dries [8] and Dan Miller [28] have worked with {\bf convergent  Weierstra\ss $ $ systems} on the reals. We consider convergent Weierstra\ss $ $ systems on subfields of the reals. For example, the algebraic power series form a convergent Weierstra\ss $ $ system on the field of real algebraic numbers. We show more generally that the coefficients of the power series in a convergent Weierstra\ss $ $ system form a real closed field, called its coefficient field.
	The system of algebraic power series does not allow formal integration.
	We introduce convergent systems with formal integration, called {\bf convergent Analysis systems}.
	We show that a convergent Analysis system is a convergent Weierstra\ss $ $ system.
	Hence given a convergent Weierstra\ss $ $ system we can ask for a smallest convergent Analysis system containing the former. In our setting this closure operation can be defined. It is important to realize that the coefficient field can change here. We show, using the approach of [19], that the system of integrated algebraic power series is exactly the closure in that  sense of the system of algebraic power series. Combining these results we obtain Theorem A.
	
	For Theorem B we are working in the setting of o-minimality and use and generalize deep results of this area as follows (compare [19] for a similar result on the reals, starting with the power series that are algebraic over the real polynomial ring). We fix a convergent Weierstra\ss $ $ system with coefficient field not necessarily the reals. A {\bf restricted power series} is given on the unit cube by a power series in this system converging in a neighbourhood thereof and extend it by zero to the whole space. Consider the structure generated by these functions expanding the coefficient field. We generalize the work of Denef and van den Dries [6] and van den Dries, Marker and Macintyre [10] to show the following. The structure is o-minimal, has quantifier elimination in the language having symbols for the restricted functions and the reciprocal and allows universal axiomatization in this language extended by symbols for power functions with rational exponents. For convergent Weierstra\ss $ $ systems on the reals this was also done by Rambaud [30] in his thesis. We have the additional difficulty that if the coefficient field is not the field of reals we have to circumvent compactness arguments.
	The results on quantifier elimination and universal axiomatization imply that definable functions are given piecewise by terms in the last language above.
	With this in hand we can extend the work of Lion and Rolin [26] to obtain preparation for definable functions. Note that in [28] Lion-Rolin preparation has been established for systems on the reals where only Weierstra\ss $ $ preparation is used. 
	If our system is a convergent Analysis system we can extend the results of Lion and Rolin [27] and  Comte, Lion and Rolin [5] on integration (see also Cluckers and Dan Miller [3] for an extension and [21, 22] for an application to integration on non-standard real closed fields) to show that parameterized integrals in the structure of restricted power series with respect to the system are given by products of sums of definable functions and logarithms of positive definable functions. This is well-defined since the coefficient field of a convergent Analysis system is closed under logarithm.
	Combining and applying these results to the convergent Analysis system of integrated algebraic series we obtain Theorem B. We can also describe families of periods in our setting. 
	
	The structure of restricted integrated algebraic functions with its good description of parameterized integrals might serve as a suitable framework for developing de Rham theory (especially if one works in the setting of constructible functions from [3] adapted to our setting). In the semialgebraic case there is an approach by Kontsevich and Soibelman [24, Appendix 8] extended and rendered by Hardt et al. [14]. There the technical setup is rather demanding since for integration the semialgebraic category is way too small.

	As mentioned above it is believed that the Euler number $e$ is not a period. But it is an integrated algebraic number since the coefficient field of a convergent Analysis system is also closed under exponentiation. We conjecture that the Euler constant $\gamma$ is not an integrated algebraic number. To capture this number we expand the structure of restricted integrated algebraic functions by the global exponential. We can extend the results of [10] to show that this structure is o-minimal, and has, in the above language extended by symbols for exponentiation and logarithm, quantifier elimination, universal axiomatization and piecewise description of definable functions by terms. This structure might also serve as a finer framework for period mappings if the varieties are defined over the rationals (compare with Bakker et al. [1]).
	
	So far we have respectively define:
	
	\vs{0.5cm}
	{\it \begin{tabular}{lcl}
			Period &=& Integral of semialgebraic function over $\IQ$\\
			&=& Integral of function definable in the field\\ 
			&& of algebraic numbers\\
			&=& Integral of function definable in the\\
			&& structure of restricted  algebraic power series \\
			&&\\
			Integrated &=& Integral of function definable in the structure\\
			algebraic number&&of restricted integrated algebraic power series\\
			&&\\
			Exponential Integrated &=& Integral of function definable in the structure\\
			algebraic number&&of restricted integrated algebraic power series \\
			&&with exponentiation\\
	\end{tabular}}
	
	\vs{0.5cm}
	On the right there are three o-minimal structures that are countable (meaning that there are only countably many definable sets (even defined with parameters)) and have excellent geometric and model theoretic properties (in particular definable functions have an explicit description).
	On the left there are subrings of the reals defined as integrals of definable functions. In the first case we have the periods. In the second case, by the above results, we can completely eliminate integration, the field of integrated algebraic numbers is precisely the universe of this structure.
	\enlargethispage{1cm}
	For the ring of {\bf exponential integrated algebraic numbers} we can show the following:
	
	\vs{0.5cm}
	{\bf Theorem C}
	
	\vs{0.1cm} 
	{\it The countable ring of exponential integrated algebraic numbers contains the exponential periods and the Euler constant $\gamma$.} 
	
	\vs{0.5cm}
	Here one has to take care what is meant by an exponential period. For  definitions of {\bf exponential periods} see [14, 25, 29] and Commelin, Habegger and Huber [4]. We follow the latter.

	In view of the above result the ring of exponential integrated algebraic numbers is a good candidate for a natural number system extending the ring of periods and containing all important mathematical constants (compare with the end of [25]).\\
	
	The paper is organized as follows. After giving some notations used through-\- out we define in Section 1 convergent power series systems over subfields of the reals as convergent Weierstra\ss $ $ systems and convergent Analysis systems and show the connections between them. In Section 2 we define integrated algebraic series and show that they form the smallest convergent Analysis system. We obtain Theorem A.
	In Section 3 we establish the model theoretic and geometric results for the structure of restricted functions from a convergent Weierstra\ss $ $ system and for the structure of restricted functions from a convergent Analysis system with exponentiation.
	In Section 4 we deal with parameterized integrals of functions definable in the structure of restricted functions from a convergent Analysis system. 
	With this we can prove Theorem B in Section 5. We introduce exponential integrated algebraic numbers and discuss the connection to exponential periods and important mathematical constants, obtaining Theorem C.
	
	\section*{Notations}
	
	By $\IN=\{1,2,3,\ldots\}$ we denote the set of natural numbers and by $\IN_0$ the set of natural numbers with $0$.
	
	Let $\IR^\infty:=\IR\cup\{\pm\infty\}$. For $R\in(\IR^\infty_{>0})^n$ we set
	$D^n(R)=\prod_{i=1}^n]-R_i,R_i[$ 
	where $]a,b[:=\{x\in \IR\mid a<x<b\}$ for $a,b\in \IR^\infty$ with $a<b$ 
	and
	$B^n(R)=\prod_{i=1}^n B(0,R_i)$
	where $B(a,r):=\{z\in\IC\mid |z-a|<r\}$ for $a\in\IC$ and $r\in \IR^\infty_{>0}$.
	
	Given a subfield $\ik$ of $\IR$ we consider always the ordering on $\ik$ coming from the reals. We set $\ik_{>0}:=\{x\in \ik\mid x>0\}$.
	
	For $n\in \IN$ we denote we denote by $\mathfrak{S}_n$ the permutation group of $\{1,\ldots,n\}$.
	
	We refer to Bochnak, Coste and Roy [2] on semialgebraic geometry and to van den Dries [9] on o-minimal structures.

	\section{Power Series Systems}
	
	\subsection{Convergent Power Series}
	
	We refer to Ruiz [31] for the theory of power series.
	
	\vs{0.5cm}
	For $n\in \IN_0$ we denote by $\ma{O}_n$ the ring of {\bf real convergent power series} in $n$ variables $X_1,\ldots, X_n$. Note that $\ma{O}_0$ is the field of reals. For $n>0$, the units of $\ma{O}_n$ are the convergent power series not vanishing in the origin. We identify an element of $\ma{O}_n$ in the natural way with an element in $\ma{O}_{n+1}$ without further notification.

	Let $f=\sum_{\alpha\in \IN_0^n}a_\alpha X^\alpha\in \ma{O}_n$ where $X=(X_1,\ldots,X_n)$ and $X^\alpha=X_1^{\alpha_1}\cdot\ldots\cdot X_n^{\alpha_n}$ for $\alpha\in \IN_0^n$.  We denote by
	$$\ma{R}(f):=\big\{R\in (\IR^\infty_{>0})^n\;\big\vert\; f\mbox{ converges (absolutely) on }D^n(R)\big\}$$
	the set of radii of convergence of $f$. 
	The set $\ma{D}(f):=\bigcup_{R\in \ma{R}(f)}D^n(R)$ is the {\bf real domain of convergence} and the set $\ma{B}(f):=\bigcup_{R\in \ma{R}(f)}B^n(R)$ is the {\bf domain of convergence} of $f$, respectively. We obtain a real analytic function $f:\ma{D}(f)\to \IR, x\mapsto \sum_{\alpha\in \IN_0^n}a_\alpha x^\alpha,$ and a holomorphic function $f_\IC:\ma{B}(f)\to \IC, z\mapsto \sum_{\alpha\in \IN_0^n}a_\alpha z^\alpha$.
	
	Let $i\in \{1,\ldots,n\}$. The formal {\bf derivative} of $f$ with respect to $X_i$ is defined by
	$\partial f/\partial X_i:=\sum_{\alpha\in \IN_0^n}\alpha_i a_\alpha X^{\alpha-\mathbf{e}_i}$ where $\mathbf{e}_i$ denotes the $i^\mathrm{th}$ unit vector.
	Given $\beta=(\beta_1,\ldots,\beta_n)\in \IN_0^n$ the formal {\bf derivative}  of order $\beta$ of $f$ is defined inductively by 
	$D^\beta(f):=\partial^{\beta_1+\ldots+\beta_n}f/\partial X_1^{\beta_1}\cdot\ldots\cdot \partial X_n^{\beta_n}$. Note that this coincides with the function germ at the origin of the derivative of order $\beta$ of the real analytic function above. Note also that
	$\ma{D}(D^\beta(f))=\ma{D}(f)$. 
	
	Let $i\in \{1,\ldots,n\}$. The formal {\bf antiderivative} of $f$ with respect to $X_i$ is defined by
	$\int f \,dX_i:=\sum_{\alpha\in \IN_0^n} \big(a_\alpha/(\alpha_i+1)\big) X^{\alpha+\mathbf{e}_i}$. Note that the antiderivative has been chosen that $(\int f\,dX_n)(X',0)=0$ (where $X':=(X_1,\ldots,X_{n-1})$), similarly for arbitrary $i$. Note also that $\ma{D}(\int f\,dX_i)=\ma{D}(f)$.

	\vs{0.5cm}
	{\bf Translation:}
	
	\vs{0.1cm}
	We need the following operation on $\ma{O}_n$.
	Let $f\in\ma{O}_n$. 
	Given $a\in \ma{D}(f)$, the germ $f_a$ of the function defined by
	$f(X+a)=f(a_1+X_1,\ldots,a_n+X_n)$ at the origin is an element of $\ma{O}_n$. We call it the {\bf translation} of $f$ at $a$ and denote it also by $f(X+a)$.
	We have
	$f_a=\sum_{\alpha\in \IN_0^n}(D^\alpha f(a)/\alpha !) X^\alpha.$
	Let $R\in \ma{R}(f)$ and let $a=(a_1,\ldots,a_n)\in D^n(R)$. Set $|a|:=(|a_1|,\ldots,|a_n|)$. Then $R-|a|\in \ma{R}(f_a)$.

	\vs{0.5cm}
	{\bf Complexification:}
	
	\vs{0.1cm}
	Let $f\in \ma{O}_n$. Let $u$ denote the real and $v$ the imaginary part of $f_{\IC}$; i.e. $f_\IC(x+iy)=u(x+iy)+iv(x+iy)$. Then the germs of $u$ and $v$ at $0\in\IR^{2n}$ are convergent power series denoted by $\Re\,f_{\IC}$ and $\Im\,f_{\IC}$, respectively. Note that $\Re\,f_\IC(X_1,0,\ldots,X_n,0)=f(X)$.
	Let $R=(R_1,\ldots,R_n)\in (\IR_{>0}^\infty)^{n}$ be a radius of convergence of $f$. Then 
	$\Re\,f_\IC$ and $\Im\,f_\IC$ have radius of convergence $(R_1/\sqrt{2},R_1/\sqrt{2},\ldots,R_n/\sqrt{2},R_n/\sqrt{2})\in (\IR_{>0}^\infty)^{2n}$.
	
	\vs{0.3cm}
	By $\ik$ we denote a subfield of the reals. Examples are given by $\IQ$, the field of rationals, and by the field of real algebraic numbers which we denote by $\IA$.
	We denote the real closure of $\ik$ in $\IR$ by $\mathrm{rc}(\ik)$. For example, $\IA=\mathrm{rc}(\IQ)$. 
	
	\vs{0.3cm}
	For $n\in\IN_0$ we denote by $\ma{O}_n(\ik):=\ik\{X\}$ the ring of convergent power series in $n$ variables with coefficients in $\ik$. 
	In the case $\ik=\IR$ we omit the field and write as above $\ma{O}_n$. 
	
	\subsection{Convergent Systems}
	
	{\bf 1.1 Definition} (Collection)
	
	\vs{0.1cm}
	By a {\bf collection (of convergent power series)} we mean a sequence $(\ma{S}_n)_{n\in \IN_0}$ such that $\ma{S}_n\subset \ma{O}_n$ for all $n\in \IN_0$.
	
	\vs{0.5cm}
	The set of all collections is partially
	ordered by setting 
	$\ma{S}\leq \ma{T}$ if $\ma{S}_n\subset \ma{T}_n$ for all $n\in \IN_0$. We say then that $\ma{S}$ is contained in $\ma{T}$.
	
	We call a collection $\ma{S}$ countable if $\ma{S}_n$ is countable for every $n\in \IN_0$.
	
	For $\sigma\in \mathfrak{S}_n$ we set $\sigma(X):=(X_{\sigma(i)},\ldots,X_{\sigma(n)})$. 
	We say that a subset $O$ of some $\ma{O}_n$ is {\bf closed under permuting} variables if $f(\sigma(X))\in O$ for every $f\in O$ and $\sigma\in \mathfrak{S}_n$.
	
	\vs{0.5cm}
	{\bf 1.2 Definition} (Presystem)
	
	\vs{0.1cm}
	We call a collection $\ma{S}=(\ma{S}_n)_{n\in \IN_0}$ a {\bf presystem (of convergent power series)}  (short: PS) if the following holds for every $n\in \IN_0$:
	\begin{itemize}
		\item[(P1)] $\ma{S}_n$ is a subring of $\ma{O}_n$.
		\item[(P2)] $\ma{S}_n=\ma{S}_{n+1}\cap \ma{O}_n$. 
		\item[(P3)] $\ma{S}_n$ is closed under permuting variables.
	\end{itemize}
	
	\newpage
	{\bf 1.3 Remark}
	
	\vs{0.1cm}
	Axiom (P3) is merely for convenience. It simplifies notations below since one can use in notations the distinguished last variable.
	
	\vs{0.5cm}
	{\bf 1.4 Definition} (System)
	
	\vs{0.1cm}
	We call a presystem $\ma{S}=(\ma{S}_n)_{n\in \IN_0}$ a {\bf system (of convergent power series)}  (short: S) if the following holds for every $n\in \IN$:
	\begin{itemize}
		\item[(S1)] $\ik:=\ma{S}_0$ is a subfield of $\IR$.
		\item[(S2)] We have $\ik[X]\subset \ma{S}_n\subset \ik\{X\}$. 
	\end{itemize}
	We call $\ik=\ik(\ma{S})$ the {\bf coefficient field} of $\ma{S}$.
	
	\vs{0.5cm}
	{\bf 1.5 Example}
	
	\vs{0.1cm}
	Let $\ik$ be a subfield of $\IR$. The collection $\ma{O}(\ik)=(\ma{O}_n(\ik))_{n\in \IN_0}$ is an uncountable system.
	
	\vs{0.1cm}
	{\bf Proof:}
	
	\vs{0.1cm}
	That $\ma{O}(\ik)$ is a system is clear. We show that $\ma{O}_1(\ik)$ is uncountable and are done. Let $\mathbf{2}^{\IN_0}$ be the set of sequences $(a_p)_{p\in \IN_0}$ with $a_p\in \{0,1\}$ for all $p\in \IN_0$.
	The map
	$$\mathbf{2}^{\IN_0}\to \ma{O}_1(\ik), (a_p)_{p\in \IN_0}\mapsto \sum_{p=0}^\infty a_pX^p,$$
	is well-defined and injective. This gives the assertion.
	\hfill$\blacksquare$

	\vs{0.5cm}
	{\bf 1.6 Definition} (Convergent System)
	
	\vs{0.1cm}
	A system $\ma{S}=(\ma{S}_n)_{n\in \IN_0}$ with coefficient field $\ik$ is called a {\bf convergent system}  (short: CS) if the following holds for all $n\in \IN_0$:
	\begin{itemize}
		\item[(C)] Let $f\in \ma{S}_n$ and $a\in \ma{D}(f)\cap \ik^n$. Then $f(X+a)\in \ma{S}_n$. 
	\end{itemize}
	
	\vs{0.2cm}
	{\bf 1.7 Example}
	
	\vs{0.1cm}
	Let $\ik$ be a subfield of $\IR$. The system $\ma{O}(\ik)$ is convergent if and only if 
	$\ik=\IR$.
	
	\vs{0.1cm}
	{\bf Proof:}
	
	\vs{0.1cm}
	That $\ma{O}$ is CS is clear (see [31, I.2]). For the other direction let $\ik\neq \IR$ and assume that $\ma{O}(\ik)$ is CS. Let $\rho\in \IR\setminus \ik$ be with $0<\rho<1$.
	Let $\sum_{p=0}^\infty a_p10^{-p}$ be the decimal expansion of $\rho$.
	We have that $f:=\sum_{p=0}^\infty a_pX^p\in \ma{O}_1(\IQ)\subset \ma{O}_1(\ik)$.
	Moreover, $1/10\in \ma{D}(f)$.
	By assumption $f_{1/10}\in \ma{O}_1(\ik)$. But  $f_{1/10}(0)=\rho\notin \ik$, contradiction.
	\hfill$\blacksquare$

	\newpage
	{\bf 1.8 Example}
	
	\vs{0.1cm}
	Let $\ik$ be a subfield of $\IR$. For $n\in \IN_0$ let $\ma{N}_n(\ik)$ be the set of power series in $\IR[[X]]$ that are algebraic over the polynomial ring $\ik[X]$ in $n$ variables. The following holds for the collection $\ma{N}(\ik):=(\ma{N}_n(\ik))_{n\in\IN_0}$:
	\begin{itemize}
		\item[(1)] $\ma{N}(\ik)$ is a convergent system with coefficient field $\mathrm{rc}(\ik)$. We have $\ma{N}(\ik)=\ma{N}(\mathrm{rc}(\ik))$.
		\item[(2)] $\ma{N}(\ik)$ is countable if and only if $\ik$ is countable. 
	\end{itemize}
	
	{\bf Proof:}
	
	\vs{0.1cm}
	(1): First we assume that $\ik$ is real closed. 
	That $\ma{N}(\ik)$ is a convergent system with coefficient field $\ik$ can be found in [2, Chapter 8]. Note that $\ma{N}_n(\ik)$ is the ring of germs of Nash functions at $0\in\IR^n$ which are defined over $\ik$; i.e. of functions that are real analytic and semialgebraic over $\ik$. Given $f\in \ma{N}_n(\ik)$ the associated real analytic function is Nash over $\ik$ at every point of its domain of convergence with components in $\ik$. Hence $f_a\in \ma{N}_n(\ik)$ for every $a\in \ma{D}(f)\cap \ik^n$.
	
	In the general case we show that $\ma{N}(\ik)=\ma{N}(\mathrm{rc}(\ik))$ and are done by above. This follows since $\mathrm{rc}(\ik)[X]$ is algebraic over $\ik[X]$ and the fact that algebraicity is transitive (see Kaplansky [23, Remark after Theorem 22 in Chapter I \S 3]).
	
	\vs{0.2cm}
	(2): 
	If $\ik$ is uncountable then clearly $\mathrm{rc}(\ik)=\ma{N}_0(\ik)$ is uncountable and therefore $\ma{N}(\ik)$ is uncountable. 
	If $\ik$ is countable then so is $\ik[X]$ and hence $\ma{N}_n(\ik)$ for every $n\in \IN_0$.
	
	\hfill$\blacksquare$
	
	\vs{0.5cm}
	In the case $\ik=\IQ$ we write $\ma{A}=(\ma{A}_n)_{n\in\IN_0}$ and speak of {\bf algebraic power series}. By above the coefficient field of $\ma{A}$ is given by the field $\IA$ of real algebraic numbers. In the case $\ik=\IR$ we write $\ma{N}=(\ma{N}_n)_{n\in \IN_0}$.
	
	\vs{0.5cm}
	The following remark is straight forward but noteworthy.
	
	\vs{0.5cm}
	{\bf 1.9 Remark}
	
	\vs{0.1cm}
	Let $\ma{S}=(\ma{S}_n)_{n\in \IN_0}$ be a convergent system with coefficient field $\ik$. 
	Let $n\in \IN$ and  $f\in \ma{S}_n$. Then $D^\beta f(a)\in \ik$ for every $a\in \ma{D}(f)\cap \ik^n$ and every $\beta\in \IN_0^n$.
	
	\vs{0.1cm}
	{\bf Proof:}
	
	\vs{0.1cm}
	By Axiom (C) we have that $f_a\in \ma{S}_n$. By Axiom (S2) we have $f_a\in \ma{O}_n(\ik)$. 
	Let $a_\beta$ be the coefficient of $f_a$ at $X^\beta$. Since $D^\beta f(a)=a_\beta\cdot \beta!$ we are done.
	\hfill$\blacksquare$
	
	\newpage
	\subsection{Weierstra\ss $ $ Systems}
	
	For the following, given $n\in \IN$, we set $X=(X_1,\ldots,X_n)$ and $X':=(X_1,\ldots,X_{n-1})$.
	
	\vs{0.5cm}
	{\bf Weierstra\ss $ $ Preparation and Division:}
	
	\vs{0.1cm}
	Recall that a convergent power series $f\in \ma{O}_n$ is called {\bf regular in $X_n$} if $f(0,X_n)\neq 0$. It is called {\bf regular in $X_n$ of order $d$}
	if 
	$$f(0)=\ldots=(\partial^{d-1}f/\partial X_n^{d-1})(0)=0\mbox{ and }(\partial^d f/\partial X_n^d)(0)\neq 0.$$
	
	If $f\neq 0$ then there is a linear transformation $A$ mapping
	$(X_1,\ldots,X_n)$ to $(X_1+c_1X_n,\ldots,X_{n-1}+c_{n-1}X_n,X_n)$
	(one can $c_1,\ldots,c_{n-1}$ choose to be rational) such that $f\circ A$ is regular in $X_n$. 
	
	Recall that $f\in\ma{O}_n$ is called a {\bf Weierstra\ss $ $ polynomial  in $X_n$ of degree $d$} if
	$$f(X)=X_n^d+a_{d-1}(X')X_n^{d-1}+\ldots + a_0(X')$$ 
	and $f$ is regular in $X_n$ of order $d$ (or, equivalently, $a_{d-1}(0)=\ldots=a_{0}(0)=0$).
	
	\vs{0.3cm}
	The {\bf Weierstra\ss $ $ preparation theorem} states:
	Let $f\in\ma{O}_n$ be regular in $X_n$ of order $d$.
	Then there is a unique Weierstra\ss $ $ polynomial $P\in \ma{O}_{n-1}[X_n]$ in $X_n$ of degree $d$  and a unique $u\in \ma{O}_n$ with $u(0)\neq 0$ such that $f=Pu$.
	
	\vs{0.3cm}
	The {\bf Weierstra\ss $ $ division theorem} states:
	Let $f\in\ma{O}_n$ be regular in $X_n$ of order $d$.
	For $g\in\ma{O}_n$ there is a unique
	$h\in\ma{O}_n$ and a unique $r\in \ma{O}_{n-1}[X_n]$ of degree smaller then $d$ such that $g=fh+r$.
	
	\vs{0.5cm}
	We define a Weierstra\ss $ $ system to be a system which allows Weierstra\ss $ $ division. 
	
	\vs{0.5cm}
	{\bf 1.10 Definition} (Weierstra\ss $ $ System)
	
	\vs{0.1cm}
	We call a convergent system
	$\ma{S}=(\ma{S}_n)_{n\in\IN_0}$ with coefficient field $\ik$
	a {\bf convergent Weierstra\ss $ $ system} (short: CWS) if the following holds for every $n\in \IN$:
	\begin{itemize}
		\item[(W)] Let $f\in\ma{S}_n$ be regular in $X_n$ of order $d$.
		For $g\in\ma{S}_n$ there is 
		$h\in\ma{S}_n$ and $r(X)\in \ma{S}_{n-1}[X_n]$ of degree smaller than $d$ such that $g=fh+r$.
	\end{itemize}
	
	\vs{0.2cm}
	{\bf 1.11 Comment}
	
	\begin{itemize}
		\item[(1)] In [7, \S 1] and [8, (1.2)] there is an additional axiom on units. This is not necessary, see Remark 1.13 below.
		\item[(2)]
		In [8, (1.2)] the definition of a convergent (Weierstra\ss) system is weaker. There translation is only demanded on some open neighbourhood of the origin whereas we request it on the whole domain of convergence.
		This is inspired on the one hand by the fact that the existing systems of interest have this property and on the other hand that this allows to define closure operations below. We need it also for covering arguments, for example in the proof of Theorem A and mainly in Section 3.
	\end{itemize}
	
	\vs{0.2cm}
	{\bf 1.12 Example}
	
	\begin{itemize}
		\item[(1)]
		The system $\ma{O}$ is a convergent Weierstra\ss $ $ system.
		\item[(2)] 
		Let $\ik$ be a subfield of $\IR$. The system $\ma{N}(\ik)$ is a convergent Weierstra\ss $ $ system.
	\end{itemize}
	
	{\bf Proof:}
	
	\vs{0.1cm}
	(1): This is clear. 
	
	\vs{0.2cm}
	(2): See [2, Chapter 8 \S 2].
	\hfill$\blacksquare$
	
	\vs{0.5cm}
	A convergent Weierstra\ss $ $ system behaves algebraically very well.
	
	\vs{0.5cm}
	{\bf 1.13 Remark}
	
	\vs{0.1cm}
	Let $\ma{S}=(\ma{S}_n)_{n\in\IN_0}$ be a convergent Weierstra\ss $ $ system.
	Let $n\in \IN$ and $f\in \ma{S}_n$ be with $f(0)\neq 0$. Then $1/f\in \ma{S}_n$.
	
	\vs{0.1cm}
	{\bf Proof:}
	
	\vs{0.1cm}
	This follows by applying (W) to $g=1$ (see [20, Remark 2.3]).
	\hfill$\blacksquare$
	
	\vs{0.5cm}
	{\bf 1.14 Remark}
	
	Let $\ma{S}=(\ma{S}_n)_{n\in\IN_0}$ be a convergent Weierstra\ss $ $ system with coefficient field $\ik$.
	For every $n\in \IN$, $\ma{S}_n$ is a regular local ring of dimension $n$ with maximal ideal generated by $X$ whose completion is $\ik[[X]]$. In particular it is Noetherian and a unique factorization domain.
	
	\vs{0.1cm}
	{\bf Proof:}
	
	\vs{0.1cm}
	See for example [2, Chapter 8 \S 2].
	\hfill$\blacksquare$
	
	\vs{0.5cm}
	A convergent Weierstra\ss $ $ system fulfills various closure properties.
	
	\vs{0.5cm}
	{\bf 1.15 Proposition}
	
	\vs{0.1cm}
	{\it Let $\ma{S}=(\ma{S}_n)_{n\in\IN_0}$ be a convergent Weierstra\ss $ $ system with coefficient field $\ik$. 
		The following holds:
		\begin{itemize}
			\item[(1)]
			{\bf Weierstra\ss $ $ Preparation:} 
			Let $n\in\IN,$ and $f\in\ma{S}_n$ be regular in $X_n$ of order $d$.
			Then there is a Weierstra\ss $ $ polynomial $P\in \ma{S}_{n-1}[X_n]$ in $X_n$ of degree $d$  and $u\in \ma{S}_n$ with $u(0)\neq 0$ such that $f=Pu$.
			\item[(2)] {\bf Composition:}
			Let $n\in \IN,
			f\in \ma{S}_n$ and $h\in \ma{S}_n^n$ be with $h(0)=0$. Then $f(h(X))\in \ma{S}_n$.
			\item[(3)] {\bf Differentiation:}
			Let $n\in \IN, f\in \ma{S}_n$ and $i\in \{1,\ldots,n\}$. Then $\partial f/\partial X_i\in \ma{S}_n$.
			\item[(4)] {\bf Inverse Function:}
			Let $n\in \IN$ and $f=(f_1,\ldots,f_n)\in \ma{S}_n^n$ be with $f(0)=0$ and $\mathrm{det}(\partial f/\partial X (0))\neq 0$. Then there is a (unique) $g\in \ma{S}_n^n$ with $g(0)=0$ such that $f(g(X))=X$.
			\item[(5)] {\bf Implicit Function:}
			Let $n,l\in \IN$ and $f=f(X,Y)=(f_1,\ldots,f_l)\in \ma{S}_{n+l}^l$ (where $X=(X_1,\ldots,X_n)$ and $Y=(Y_1,\ldots,Y_l)$) be with $f(0)=0$ and 
			$\mathrm{det}(\partial f/\partial Y(0))\neq 0$.
			Then there is a (unique) $g\in \ma{S}_n^l$ with $g(0)=0$ such that $f(X,g(X))=0$.
			\item[(6)] {\bf Complexification:}
			Let $n\in \IN$ and $f\in \ma{S}_n$. Then $\Re\,f_\IC,\Im\,f_\IC\in \ma{S}_{2n}$.
	\end{itemize}}
	{\bf Proof:}
	
	\vs{0.1cm}
	(1): That the Weierstra\ss $ $ division theorem implies the Weierstra\ss $ $ preparation theorem is classical (see  for example [31, Proposition I.3.3]).
	
	\vs{0.2cm}
	(2) - (5): See [7, Remarks 1.3].
	
	\vs{0.2cm}
	(6): See [28, Proposition 9.1(i)].
	\hfill$\blacksquare$

	\vs{0.5cm}
	{\bf 1.16 Theorem}
	
	\vs{0.1cm}
	{\it The coefficient field of a convergent Weierstra\ss $ $ system is real closed.}
	
	\vs{0.1cm}
	{\bf Proof:}
	
	\vs{0.1cm}
	Denote by $\ma{S}=(\ma{S}_n)_{n\in\IN_0}$ the convergent Weierstra\ss $ $ system and by $\ik$ its coefficient field. In the proof we heavily use that $\ik$ is dense in the reals.
	Let $f\in \ik[X]$ be univariate and let $a\in \IR$ be a zero of $f$. We show that $a\in \ik$.
	This gives the proof.
	We may assume that $f'(a)\neq 0$, otherwise we consider a sufficiently high derivative of $f$.
	Let $b,c \in \ik$ be such that $b<a<c$ and $f'(x)\neq 0$ for all $x\in [b,c]$.
	Without restriction we may assume that $f$ is strictly increasing on $]b,c[$. 
	Let $d:=f(b), e:=f(c)$ and $g:[d,e]\to [b,c]$ be the inverse of $f|_{[a,b]}$. Note that $d<0<e$.
	
	\vs{0.2cm}
	{\bf Claim:} For every $y\in [d,e]$ we have a uniquely determined $h_{(y)}\in \ma{O}_1$ such that $g$ coincides with $h_{(y)}(Y-y)$ on 
	$(\ma{D}(h_{(y)})+y)\cap [d,e]$. If $y\in f([b,c]\cap \ik)$ we have $h_{(y)}\in \ma{S}_1$.
	
	\vs{0.1cm}
	{\bf Proof of the claim:}
	Uniqueness being clear we show existence.
	Let $x:=g^{-1}(y)\in [b,c]$ and let $F:=f(X+x)-y\in \IR[X]\subset \ma{O}_1$. 
	Then $F(0)=0$ and $F'(0)=f'(x)\neq 0$. By Proposition 1.15(4) applied to $\ma{O}$ we have that $F$ has a compositional inverse $G\in \ma{O}_1$ with $G(0)=0$.
	We take $h_{(y)}:=x+G$.
	If $y\in f([b,c]\cap\ik)$ then $x,y\in \ik$ and therefore $F\in \ik[X]\subset \ma{S}_1$. Applying Proposition 1.15(4) to $\ma{S}$ we obtain that $G$ and therefore $h_{(y)}$ belong to $\ma{S}_1$.
	
	\hfill$\blacksquare_\mathrm{Claim}$
	
	\newpage
	Since $f([b,c]\cap \ik)$ is dense in $[d,e]$ we obtain that the sets
	$y+\ma{D}(h_{(y)})$ with $y\in f([b,c]\cap \ik)$ cover $[b,c]$. (Note that we use here the full strength of our Axiom (C).)
	Hence we can take $\tilde{y}\in f([b,c]\cap \ik)$ such that $0=\tilde{y}+t$ for some $t\in \ma{D}(h_{(\tilde{y})})$. Then $t\in \ik$ and we obtain by Remark 1.9 that $a=g(0)=h_{(\tilde{y})}(t)\in \ik$. 
	\hfill$\blacksquare$
	
	\vs{0.5cm}
	{\bf 1.17 Corollary}
	
	\vs{0.1cm}
	{\it The system $\ma{A}$ of algebraic power series is the smallest convergent Weierstra\ss $ $ system.}
	
	\vs{0.1cm}
	{\bf Proof:}
	
	\vs{0.1cm}
	Let $\ma{S}=(\ma{S}_n)_{n\in \IN_0}$ be a convergent Weierstra\ss $ $ system and let $\ik$ be its coefficient field.
	By Theorem 1.16 the field $\ik$ is real closed. Hence it contains the field $\IA$ of real algebraic numbers.
	Let $n\in \IN$ and $f\in \ma{A}_n$. We have to show that $f\in \ma{S}_n$.
	By [2, Theorem 8.1.4] we find $l\in \IN_0$ such that $f(X)=p(X,g(X))$ where $p(X,Y)\in \IA[X,Y]$ (with $X=(X_1,\ldots,X_n)$ and $Y=(Y_1,\ldots,Y_l)$) and $g$ is implicitely defined by $g(0)=0$ and $q(X,g(X))=0$ for some $q(X,Y)\in \IA[X,Y]^l$ such that $q(0)=0$ and $\det(\partial q/\partial Y)(0)\neq 0$.
	By Axiom (S2) and Proposition 1.15(2,5) we obtain that $f\in \ma{S}_n$.
	\hfill$\blacksquare$
	
	\subsection{Analysis Systems}
	
	{\bf 1.18 Remark}
	
	\vs{0.1cm}
	A convergent Weierstra\ss $ $ system is not necessarily closed under taking formal antiderivatives.
	
	\vs{0.1cm}
	{\bf Proof:}
	
	\vs{0.1cm}
	We have that $1/(1+X)=\sum_{p=0}^\infty (-1)^pX^p\in \ma{A}_1$. Its formal antiderivative is the logarithmic series $L(X)$ which is transcendental (i.e. not algebraic over the real polynomial ring $\IR[X]$) and therefore not an algebraic power series.
	\hfill$\blacksquare$
	
	\vs{0.5cm}
	We add formal integration in the following way.
	
	\vs{0.5cm}
	{\bf 1.19 Definition} (Analysis System)
	
	\vs{0.1cm}
	We call a convergent system $\ma{S}=(\ma{S}_n)_{n\in\IN_0}$ a {\bf  convergent Analysis system} (short: CAS) if the following holds for all $n\in \IN$ and $f\in \ma{S}_n$:
	\begin{itemize}
		\item[(A1)] $1/f\in \ma{S}_n$ if $f(0)\neq 0$,
		\item[(A2)] $f(X',0)\in \ma{S}_{n-1}$, 
		\item[(A3)] $\partial f/\partial X_n\in \ma{S}_n$,
		\item[(A4)] $\Re\,f_\IC, \Im\,f_\IC\in \ma{S}_{2n}$,
		\item[(A5)] $\int f\, dX_n\in \ma{S}_n$ if $n\geq 1$.
	\end{itemize}
	
	\vs{0.2cm}
	{\bf 1.20 Remark}
	
	\vs{0.1cm}
	Let $\ma{S}=(\ma{S}_n)_{n\in \IN_0}$ be a convergent Analysis system. Let $n\in \IN$ and $f\in \ma{S}_n$. Then $\partial f/\partial X_i,\int f\,dX_i\in \ma{S}_n$ for all $i\in \{1,\ldots,n\}$.
	
	\vs{0.1cm}
	{\bf Proof:}
	
	\vs{0.1cm}
	This follows with Axiom (P3).
	\hfill$\blacksquare$
	
	\vs{0.5cm}
	In [19] we have constructed the system $\ma{IN}$ of integrated algebraic power series over the reals. 
	
	\vs{0.5cm}
	{\bf 1.21 Examples}
	
	\begin{itemize}
		\item[(1)] The system $\ma{O}$ is a CAS.
		\item[(2)]
		Let $\ik$ be a subfield of $\IR$. The system $\ma{N}(\ik)$ is not a CAS.
		\item[(3)] 
		The system $\ma{IN}$ is a CAS. 
	\end{itemize}
	{\bf Proof:}
	
	\vs{0.1cm}
	(1): This is clear.
	
	\vs{0.1cm}
	(2): See the proof of Remark 1.18.
	
	\vs{0.1cm}
	(3): This follows by [19, Def. 1.2, Prop. 2.2 \& Th. 2.7].
	\hfill$\blacksquare$
	
	\vs{0.5cm}
	We have seen that the system $\ma{A}$ of algebraic series is a convergent Weierstra\ss $ $ system with coefficient field $\IA$. 
	The coefficient field of a convergent Analysis system necessarily contains the transcendental number $\pi$. 
	
	\vs{0.5cm}
	{\bf 1.22 Proposition}
	
	\vs{0.1cm}
	{\it Let $\ma{S}=(\ma{S}_n)_{n\in\IN_0}$ be a convergent Analysis system with 
		coefficient field $\ik$.
		Then $\pi\in \ik$.}
	
	\vs{0.1cm}
	{\bf Proof:}
	
	\vs{0.1cm}
	We have $f:=1+X^2\in \ma{S}_1$ by Axiom (S2).
	By Axiom (A1) we get that $g:=1/f=\sum_{p=0}^\infty (-1)^pX^{2p}\in \ma{S}_1$. 
	By Axiom (A5) we get that $h:=\int g\,dX=\sum_{p=0}^\infty \big((-1)^p/(2p+1)\big)X^{2p+1}\in \ma{S}_1$. 
	The series $h$ is the arctangent series and has radius of convergence $1$. Since $\ik$ contains the rationals we are done by Remark 1.9 and the following formula of John Machin:
	$$\frac{\pi}{4}=4\arctan\Big(\frac{1}{5}\Big)-\arctan\Big(\frac{1}{239}\Big).$$
	\hfill$\blacksquare$
	
	\vs{0.5cm}
	In a convergent Analysis system we can perform integration in the following sense. 
	
	\newpage
	{\bf 1.23 Proposition} (Integration)
	
	\vs{0.1cm}
	{\it Let $\ma{S}=(\ma{S}_n)_{n\in\IN_0}$ be a convergent Analysis system with coefficient field $\ik$. Let $n\in \IN$ and $f\in\ma{S}_n$. 
		Let $R=(R_1,\ldots,R_n)\in \ma{R}(f)$ and let $R':=(R_1,\ldots,R_{n-1})$. We set $D:=D^n(R)$ and $D':=D^{n-1}(R')$.
		Let $a\in \ik_{>0}$ be such that $(0,a)\in D$. Let
		$$g:D'\to \IR, g(x')=\int_0^af(x',\zeta)\,d\zeta.$$
		Then the germ of $g$ at the origin belongs to $\ma{S}_{n-1}$.}
	
	\vs{0.1cm}
	{\bf Proof:}
	
	\vs{0.1cm}
	We denote the germ of $g$ again by $g$. We have 
	$g(X')=(\int f\,dX_n)_{(0,a)}(X',0)$ and are done by Axioms (C), (A2) and (A5).
	\hfill$\blacksquare$
	
	\vs{0.5cm}
	We are able to show that CAS implies CWS.
	
	\vs{0.5cm}
	{\bf 1.24 Theorem}
	
	\vs{0.1cm}
	{\it A convergent Analysis system is a convergent Weierstra\ss $ $ system.}

	\vs{0.1cm}
	{\bf Proof:}
	
	\vs{0.1cm}
	Using the previous results we can as in [19] adjust the classical proofs using complex integration (cf. for example Gunning and Rossi [12, II.B]). We give the proof that the reader sees how the axioms are involved and where $\pi$ is needed.
	Let $\ma{S}=(\ma{S}_n)_{n\in \IN_0}$ be a convergent Analysis system.
	
	\vs{0.5cm}
	{\bf Step 1: Weierstra\ss $ $ Preparation}
	
	\vs{0.1cm}
	We first establish Weierstra\ss $ $ preparation for $\ma{S}$. Let $n\in \IN$ and let $f\in \ma{S}_n$ be regular in $X_n$ of order $d$.
	Let $R=(R_1,\ldots,R_n)$ be a radius of convergence of $f$. We do not distinguish in the following between $f$ and $f_\IC$. We indicate the difference by writing $f(x)$ respectively $f(z)$. We also use frequently the identification of $\IC$ with $\IR^2$.
	Since $f$ is regular of order $d$ we find $\delta=(\delta_1,\ldots,\delta_n)\in \ik_{>0}^n$ such that $\sqrt{2}\delta_i<R_i$ for all $i\in\{1,\ldots,n\}$ and the following holds:
	For every $z'\in B^{n-1}(\delta')$ where $\delta':=(\delta_1,\ldots,\delta_{n-1})\in \ik_{>0}^{n-1}$ the function $f_{(z')}:B(0,R_n)\to\IC, z_n\mapsto f(z',z_n)$, has exactly $d$ zeros (counted with multiplicities) and these are contained in $B(0,\delta_n)$. The zeroes will be denoted by $\varphi_1(z'),\ldots,\varphi_d(z').$ Let
	$$
	P:B^{n-1}(\delta')\times \IC\to\IC, P(z',z_n)=\prod_{r=1}^d\left(z_n-\varphi_r(z')\right).
	$$
	Then by the classical proof $P(x)\in\ma{O}_{n-1}[X_n]$ is the unique Weierstra\ss $ $ polynomial for $f$. We show that $P\in \ma{S}_{n-1}[X_n]$.
	On $B^{n-1}(\delta')\times\IC$ the function $P$ can be written as $P(z)=z_n^d+a_1(z')z_n^{d-1}+\ldots+a_d(z')$ where $a_k(z')$ are the elementary symmetric functions in the values $\varphi_r(z')$. The functions $a_k(z')$ are polynomials (with rational coefficients) in the power sums  
	$$h_s:B^{n-1}(\delta')\to \IC, z'\mapsto \sum_{r=1}^d\varphi_r(z')^s,$$
	where $s\in \{1,\ldots,d\}$. So it is enough to show that $h_1(x),\ldots,h_d(x)\in\ma{S}_{n-1}$. We fix $s\in \{1,\ldots,d\}$. One has the integral representation
	$$h_s(z')=\frac{1}{2\pi i}\int_{\partial A}\frac{\partial f(z',\zeta)}{\partial\zeta}\frac{\zeta^s}{f(z',\zeta)}d\zeta$$
	where $A:=[-\delta_n,\delta_n]^2$.
	Note that $B(0,\delta_n)\subset \mathring{A}\subset B(0,R_n)$ and that $f(z',\zeta)\neq 0$ for $z'\in B^{n-1}(\delta')$ and
	$\zeta\in \partial A$. 
	Let
	$$\Psi_s:B^{n-1}(\delta')\times \partial A\to \IC, (z',\zeta)\mapsto \frac{\partial f(z',\zeta)}{\partial \zeta}\frac{\zeta^s}{f(z',\zeta)},$$
	be the integrand. 
	Using Axioms (A1), (A3) and (A4) we find for every $\eta\in \partial A\cap \ik^2$ power series $F^{(\eta)},G^{(\eta)}\in \ma{S}_{2n}$ such that the following holds 
	\begin{itemize}
		\item[(1)]  $\Psi_s(z',\zeta)=F^{(\eta)}(z',\zeta-\eta)+i\, G^{(\eta)}(z',\zeta-\eta)$ for all $(z',\zeta)\in D^{(\eta)}$ with $\zeta\in \partial A$.
	\end{itemize}
	where 
	$$D^{(\eta)}:=\big(\ma{D}(F^{(\eta)})\cap \ma{D}(G^{(\eta)})\big)+(0,\eta).$$ 
	Note that on the right side we view $(z',\zeta-\eta)$ as an element of $\IR^{2n}$.
	Since $\delta_n\in \ik$ we have that $\ik^2$ is dense in $\partial A$. Using that $\partial A$ is compact and the full strength of Axiom (C) (compare with the proof of Theorem 1.16) we obtain the existence of a finite subset $\ma{E}$ of $\partial A$ contained in $\ik^2$ and for every $\eta\in \ma{E}$ some horizontal or vertical closed segment $I^{(\eta)}$ contained in $\partial A$ such that the following holds:
	\begin{itemize}
		\item[(2)] There is some $r^{(\eta)}\in \IR_{>0}^n$ such that $I^{(\eta)}\subset B^n(r^{(\eta)})\subset D^{(\eta)}$, 
		\item[(3)] the endpoints of $I^{(\eta)}$ are contained in $\ik^2$ and one of them is $\eta$,
		\item[(4)] $I^{(\eta)}$ and $I^{(\eta')}$ have at most one point in common for distinct $\eta,\eta'\in \ma{E}$,
		\item[(5)] $\bigcup_{\eta\in\ma{E}}I^{(\eta)}=\partial A$.
	\end{itemize}
	Let $r'\in \IR_{>0}^{n-1}$ be such that 
	$$B':=B^{n-1}(r')\subset \bigcap_{\eta\in \ma{E}} B^{n-1}((r^{(\eta)})')$$
	where $(r^{(\eta)})':=(r_1^{(\eta)},\ldots,r_{n-1}^{(\eta)})$ for $\eta\in \ma{E}$.
	We obtain for $z'\in B'$ that
	$h_s(z')=\sum_{\eta\in \ma{E}}\ma{I}^{(\eta)}(z)$
	where
	$$\ma{I}^{(\eta)}:B'\to \IC, z'\mapsto \int_{I^{(\eta)}} \big(F^{(\eta)}(z',\zeta-\eta)+i\, G^{(\eta)}(z',\zeta-\eta)\big)\,d\zeta$$
	for $\eta\in \ma{E}$.
	Fix such an $\eta$ and let $c^{(\eta)}$ be the other endpoint of $I^{(\eta)}$. Setting $a^{(\eta)}:=|c^{(\eta)}-\eta|\in \ik_{\geq 0}$ and choosing the linear parameterization $[0,a^{(\eta)}]\to I^{(\eta)}$ mapping $0$ to $\eta$ we find by Axiom (C) power series $\hat{F}^{(\eta)},\hat{G}^{(\eta)}\in \ma{S}_{2n}$ with $r^{(\eta)}\in \ma{R}(\hat{F}^{(\eta)})\cap \ma{R}(\hat{G}^{(\eta)})$ such that
	$$\ma{I}^{(\eta)}(z')=\int_0^{a^{(\eta)}} \big(\hat{F}^{(\eta)}(z',\zeta)+i\,\hat{G}^{(\eta)}(z',\zeta)\big)\,d\zeta$$
	for $z'\in B'$.
	By Proposition 1.23 we get that $\ma{I}^{(\eta)}(z')\in \ma{S}_{2(n-1)}\oplus i\,\ma{S}_{2(n-1)}$. 
	This gives that $h_s(z')\in \ma{S}_{2(n-1)}\oplus i\,\ma{S}_{2(n-1)}$.
	But since $h_s(x')\in \ma{O}_{n-1}$ we obtain that 
	$h_s(x')\in \ma{S}_{n-1}$.
	
	We see by the classical proof that $u:=f/P$ is holomorphic on $B^{n-1}(\delta')\times B(0, R_n)$ and non-vanishing. It remains to show that $u\in\ma{S}_n$.
	By construction $P(z',\zeta)\neq 0$ for all $z'\in B^{n-1}(\delta')$ and all $\zeta\in \partial [-\delta_n,\delta_n]^2$. By the Cauchy integral formula we get
	$$u(z)=\frac{1}{2\pi i}\int_{\partial [-\delta_n,\delta_n]^2}\frac{f(z',\zeta)/P(z',\zeta)}{\zeta -z_n}d\zeta$$
	for all $z\in B^n(\delta)$. As above we get that $u\in\ma{S}_n$.
	
	\vs{0.5cm}
	{\bf Step 2: Weierstra\ss $ $ Division}
	
	\vs{0.1cm}
	By the preparation theorem and Axiom (A1) we can assume that $f$ is a Weierstra\ss $ $ polynomial.
	Let $R=(R_1,\ldots,R_n)$ be a radius of convergence of both $f$ and $g$ and let $\delta\in \ik_{>0}^n$ be such that $\sqrt{2}\delta_i<R_i$  for all $j\in \{1,\ldots,n\}$ and the following holds: $f(z',\zeta)\neq 0$ for all $z'\in B^{n-1}(\delta')$ and all $\zeta\in \partial [-\delta_n,\delta_n]^2$. The function
	$$h:B^n(\delta)\to\IC, h(z)=\frac{1}{2\pi i}\int_{\partial [-\delta_n,\delta_n]}\frac{g(z',\zeta)/f(z',\zeta)}{\zeta-z_n}d\zeta,$$
	is holomorphic. The function $r:=g-fh:B^n(\delta)\to\IC$ is also holomorphic and has the integral expansion
	$$r(z)=\frac{1}{2\pi i}\int_{\partial [-\delta_n,\delta_n]}
	\frac{g(z',\zeta)}{f(z',\zeta)}\left(\frac{f(z',\zeta)-f(z',z_n)}{\zeta-z_n}\right)d\zeta.$$
	As above we see that $h, r\in\ma{S}_n$. By the classical proof $h,r$ fulfill the requirements.
	\hfill$\blacksquare$

	\subsection{Closures}

	{\bf 1.25 Definition}
	
	\vs{0.1cm}
	Let $\ma{S}=(\ma{S}_n)_{n\in \IN_0}$ be a collection.
	
	\begin{itemize}
		\item[(a)] We call the smallest presystem containing $\ma{S}$ the PS-closure of $\ma{S}$ and write $\mathrm{PS}(\ma{S})$.
		\item[(b)] We call the smallest system containing $\ma{S}$ the S-closure of $\ma{S}$ and write $\mathrm{S}(\ma{S})$.
		\item[(c)] We call the smallest convergent system containing $\ma{S}$ the CS-closure of $\ma{S}$ and write $\mathrm{CS}(\ma{S})$. 
		\item[(d)] We call the smallest convergent Weierstra\ss $ $ system  containing $\ma{S}$ the CWS-closure of $\ma{S}$ and write $\mathrm{CWS}(\ma{S})$.
		\item[(e)] We call the smallest convergent Analysis system containing $\ma{S}$ the CAS-closure of $\ma{S}$ and write $\mathrm{CAS}(\ma{S})$.  
	\end{itemize}
	
	\vs{0.2cm}
	{\bf 1.26 Remark}
	
	\vs{0.1cm}
	Let $\ma{S}$ be a collection.
	\begin{itemize}
		\item[(1)] 
		The above closures exist.
		\item[(2)]
		We have 
		$$\ma{S}\leq \mathrm{PS}(\ma{S})\leq \mathrm{S}(\ma{S})\leq \mathrm{CS}(\ma{S})\leq\mathrm{CWS}(\ma{S})\leq\mathrm{CAS}(\ma{S}).$$
	\end{itemize}
	
	\vs{0.2cm}
	{\bf 1.27 Example}
	
	\vs{0.2cm}
	We have $\mathrm{CAS}(\ma{N})=\ma{IN}$.
	
	\vs{0.1cm}
	{\bf Proof:}
	
	\vs{0.1cm}
	By Example 1.21(3) we know that $\ma{IN}$ is CAS.
	By the definition of $\ma{IN}$ (see [19]) clearly $\ma{IN}$ is contained in $\mathrm{CAS}(\ma{N})$.
	\hfill$\blacksquare$

	\vs{0.5cm}
	{\bf 1.28 Proposition}
	
	\vs{0.1cm}
	{\it Let $\ma{S}$ be a countable collection. Then $\mathrm{CAS}(\ma{S})$ is countable.}
	
	\vs{0.1cm}
	{\bf Proof:}
	
	\vs{0.1cm}
	This follows from the following elementary set-theoretic argument:
	Let $X$ be a set and let $\ma{F}$ be a countable set of finitary operations on $X$ (that is, for each $F\in \ma{F}$ there is $k\in \IN$ such that $F$ is a map from a subset of $X^k$ into $X$), then the $\ma{F}$-closure in $X$ of any countable subset of $X$ is again countable.
	
	Let $\ma{T}:=\mathrm{CAS}(\ma{S})$. We let $X$ be the disjoint union of the $\ma{O}_n$ and $A$ be the disjoint union of the $\ma{S}_n$. Then $A$ is countable by assumption. There is a countable set $\ma{F}$ of finitary operations on $X$ such that the union of the $\ma{T}_n$ is the $\ma{F}$-closure of $A$ and we are done.
	\hfill$\blacksquare$

	\section{Integrated Algebraic Series and Integrated Algebraic Numbers}
	
	\subsection{Definition of Integrated Algebraic Series}
	
	We define the relevant class of power series which we call {\bf integrated algebraic power series}.
	
	\vs{0.5cm}
	{\bf 2.1 Definition}
	
	\vs{0.1cm}
	We denote by $\ma{IA}=(\ma{IA})_{n\in\IN_0}$  the smallest presystem such that the following properties hold for all $n\in\IN_0$ and $f\in \ma{IA}_n$:
	\begin{itemize}
		\item[(I)] $\ma{A}_n\subset\ma{IA}_n$.
		\item[(II)] $1/f\in\ma{IA}_n$ if $f(0)\neq 0$.
		\item[(III)] $f(h(X))\in \ma{IA}_n$ for all $h\in \IQ[X]^n$ with $h(0)=0$.
		\item[(IV)] $f(X+a)\in\ma{IA}_n$ for all $a\in \ma{D}_f\cap \IQ^n$.
		\item[(V)] $\int f\,dX_n\in \ma{IA}_n$ if $n\geq 1$.
	\end{itemize}
	
	\vs{0.2cm}
	Our first main result is the following.
	
	\vs{0.5cm}
	{\bf 2.2 Theorem}
	
	\vs{0.1cm}
	{\it The presystem $\ma{IA}$ of integrated algebraic series is the CAS-closure of the convergent Weierstra\ss $ $ system $\ma{A}$ of algebraic power series.}
	
	\vs{0.5cm}
	One part of Theorem 2.2 is evident.
	
	\vs{0.5cm}
	{\bf 2.3 Remark}
	
	\vs{0.1cm}
	$\ma{IA}$ is contained in $\mathrm{CAS}(\ma{A})$.
	
	\vs{0.1cm}
	{\bf Proof:}
	
	\vs{0.1cm}
	This follows from Theorem 1.24, Remark 1.13 and Proposition 1.15(2).
	\hfill$\blacksquare$
	
	\vs{0.5cm}
	For the other part we
	need some preparation.

	\vs{0.5cm}
	{\bf 2.4 Proposition}
	
	\vs{0.1cm}
	{\it We have that $\ik:=\ma{IA}_0$ is a field with $\pi\in \ik$.}
	
	\vs{0.1cm}
	{\bf Proof:}
	
	\vs{0.1cm}
	That $\ik$ is a field follows from (P1) and Property (II).
	To obtain that $\pi\in \ik$ we can use the proof of Proposition 1.22. Note that Properties (I) -- (V) are sufficient.
	\hfill$\blacksquare$
	
	\vs{0.5cm}
	For the rest of this section we denote by $\ik$ this field.
	
	\vs{0.5cm}
	To show that a certain property $(*)$ holds for all elements of $\ma{IA}_n$ and all $n\in\IN_0$ it is enough by the defining Properties (I) - (V) to show the following steps for all $n\in \IN_0$. 
	\begin{itemize}
		\item[($\mathrm{I}_*$)] If $f\in \ma{A}_n$ then $f$ has property $(*)$.
		\item[($\mathrm{II}_*$)] If $f,g\in\ma{IA}_n$ have property $(*)$ then $f+g, fg$ and, for $f(0)\neq 0$, $1/f$ have property $(*)$. If $f\in \ma{IA}_n$ and $\sigma\in \mathfrak{S}_n$ then $f(\sigma(X))$ has property $(*)$.
		\item[($\mathrm{III}_*$)] If $f\in\ma{IA}_n$ has property $(*)$ then $f(h(X))$ has property $(*)$ for every $h\in \IQ[X]^n$ with $h(0)=0$.
		\item[($\mathrm{IV}_*$)] If $f\in\ma{IA}_n$ has property $(*)$ then $f_a$ has property $(*)$
		for all $a\in \ma{D}(f)\cap \IQ^n$.
		\item[($\mathrm{V}_*$)] If $f\in\ma{IA}_n$ has property $(*)$ then $\int f\,dX_n$ has property $(*)$.
	\end{itemize}
	
	\vs{0.2cm}
	{\bf 2.5 Proposition}
	
	\vs{0.1cm}
	{\it $\ma{IA}_n$ is closed under differentiation.}
	
	\vs{0.1cm}
	{\bf Proof:}
	
	\vs{0.1cm}
	We can adopt the proof of [19, Proposition 2.2] as follows.
	We show the property
	$$(\mathrm{D}): f\in\ma{IA}_n\Longrightarrow \frac{\partial f}{\partial X_i}\in\ma{IA}_n\mbox{ for all }i\in \{1,\ldots,n\}.$$
	((D) stands for 'differentiation'.)
	
	\vs{0.2cm}
	$(\mathrm{I})_{\mathrm{D}}$: This is clear since $\ma{A}_n$ is closed under differentiation (see [2, Chapter II \S 9]). 
	
	\vs{0.2cm}
	For the following we assume that $f,g\in \ma{IA}_n$ have the property (D).
	
	\vs{0.2cm}
	$(\mathrm{II})_{\mathrm{D}}$: Let $i\in \{1,\ldots,n\}$. For the sum, the product and the reciprocal this follows from the rules for differentiation: We have
	$$\partial(f+g)/\partial X_i=\partial f/\partial X_i+\partial g/\partial X_i,$$ 
	$$\partial (fg)/\partial X_i=(\partial f/\partial X_i)g+f(\partial g/\partial X_i)$$ 
	$$\partial(1/f)/\partial X_i=-(\partial f/\partial X_i)/f^2$$ 
	where in the third formula we assume that $f(0)\neq 0$.
	For permutations it follows from 
	$$\partial f(\sigma(X))/\partial X_i=\partial f/\partial X_{\sigma(i)}(\sigma(X))$$ 
	for $\sigma\in\mathfrak{S}_n$.
	
	\vs{0.2cm}
	$(\mathrm{III})_{\mathrm{D}}$:
	Let $h=(h_1,\ldots,h_n)\in \IQ[X]^n$ be with $h(0)=0$. Let $i\in \{1,\ldots,n\}$.
	We have by the chain rule
	$$\frac{\partial (f(h(X)))}{\partial X_i}=\sum_{j=1}^n\left(\frac{\partial f}{\partial X_j}(h(X))\right)\frac{\partial h_j}{\partial X_i}.$$
	This shows $(\mathrm{III})_{\mathrm{D}}$.
	
	\vs{0.2cm}
	$(\mathrm{IV})_{\mathrm{D}}$:  Let $i\in \{1,\ldots,n\}$. We have $\ma{D}(\partial f/\partial X_i)=\ma{D}(f)$. 
	For $a\in \ma{D}(f)\cap \IQ^n$ we have
	$\partial f_a/\partial X_i=(\partial f/\partial X_i)_a$. Hence $\partial f_a/\partial X_i\in\ma{IA}_n$.
	
	\vs{0.2cm}
	$(\mathrm{V})_{\mathrm{D}}$: Let $i,j\in\{1,\ldots,n\}$.
	We have
	$$\frac{\partial(\int f\,dX_j)}{\partial X_i}=\left\{\begin{array}{ccc}
		f& &i=j,\\
		& \mbox{ if }&\\
		\int (\partial f/\partial X_i)\,dX_j& &i\neq j.
	\end{array}\right.
	$$
	This shows $(\mathrm{V})_{\mathrm{D}}$.
	\hfill$\blacksquare$
	
	\vs{0.5cm}
	{\bf 2.6 Corollary}
	
	\vs{0.1cm}
	{\it $\ma{IA}$ is a system with coefficient field $\ik$.}
	
	\vs{0.1cm}
	{\bf Proof:}
	
	\vs{0.1cm}
	By Proposition 2.4 $\ik$ is a field.
	Let $n\in \IN$.
	Since $\IZ[X]\subset \ma{A}_n\subset \ma{IA}_n$ we obtain by (P1) and (P2) that $\ik[X]\subset \ma{IA}_n$.
	It remains to show that $\ma{IA}_n\subset \ik\{X\}$.
	Let $f=\sum_{\alpha\in \IN_0^n}a_\alpha X^\alpha\in \ma{IA}_n$.
	For $\beta\in \IN_0^n$ we have to show that $a_\beta\in \ik$.
	We have
	$a_\beta=(D^\beta f)(0)/\beta !$ and are done by Proposition 2.5 and Property (III).
	\hfill$\blacksquare$
	
	\vs{0.5cm}
	For the following we set $Z:=(X_1,Y_1,\ldots,X_n,Y_n)$ and $Z':=(X_1,Y_1,\ldots,X_{n-1},Y_{n-1})$.
	
	\vs{0.5cm}
	{\bf 2.7 Proposition}
	
	\vs{0.1cm}
	{\it $\ma{IA}$ is closed under complexification}
	
	\vs{0.1cm}
	{\bf Proof:}
	
	\vs{0.1cm}
	We can adjust the proof of [19, Theorem A]. 
	We have to show for every $n\in \IN$ the following property.
	\begin{eqnarray*}
		(\IC):f\in\ma{IA}_n&\Longrightarrow& \Re\,f_\IC, \Im\,f_\IC\in\ma{IA}_{2n}.
	\end{eqnarray*}
	(($\IC$) stands for 'complexification'.)
	
	\vs{0.2cm}
	$(\mathrm{I}_{\IC})$:
	This follows from Example 1.12(2) and Proposition 1.15(6).
	
	\vs{0.2cm}
	For the following we assume that $f,g\in \ma{IA}_n$ have the property ($\IC$).
	
	\vs{0.2cm}
	$(\mathrm{II}_{\IC})$: It is clear that
	$$\Re(f+g)_\IC=\Re\,f_\IC+\Re\,g_\IC, \;\Im(f+g)_\IC=\Im\,f_\IC+\Im\,g_\IC$$
	and that
	$$\Re(fg)_\IC=\Re\,f_\IC \Re\,g_\IC-\Im\,f_\IC\Im\,g_\IC,$$
	$$\Im(fg)_\IC=\Re\,f_\IC\Im\,g_\IC+\Im\,f_\IC\Re\,g_\IC.$$
	This shows that $\Re(f+g)_\IC, \Im(f+g)_\IC, \Re(fg)_\IC, \Im(fg)_\IC\in\ma{IA}_{2n}$.
	
	Let $f(0)\neq 0$. Then $\Re\,f_\IC(0)\neq 0$. We have
	$$\Re\,(1/f)_\IC=\frac{\Re\,f_\IC}{(\Re\,f_\IC)^2+(\Im\,f_\IC)^2}, \;\Im\,(1/f)_\IC=-\frac{\Im\,f_\IC}{(\Re\,f_\IC)^2+(\Im\,f_\IC)^2}.$$
	So $\Re\,(1/f)_\IC, \Im\,(1/f)_\IC\in\ma{IA}_{2n}$.
	
	Let $\sigma\in \mathfrak{S}_n$. Define $\tau:\{1,\ldots,2n\}\to\{1,\ldots,2n\}$ by $\tau(i)=2\sigma((i+1)/2)-1$ for $i\in\{1,\ldots,2n\}$ odd and $\tau(i)=2\sigma(i/2)$ for $i\in \{1,\ldots,2n\}$ even. Then $\tau\in \mathfrak{S}_{2n}$.
	We have
	$\Re\,f(\sigma(X))_\IC=\Re\,f_\IC(\tau(Z))$ and $\Im\,f(\sigma(X))_\IC=\Im\,f_\IC(\tau(Z))$. Hence $\Re\,f(\sigma(X))_\IC, \Im\,f(\sigma(X))_\IC\in \ma{IA}_{2n}$.
	
	\vs{0.2cm}
	$(\mathrm{III}_{\IC})$:
	Let $h\in\IQ[X]^n$ be with $h(0)=0$. Then $\Re\,h_\IC,\Im\,h_\IC\in\IQ[Z]^n$ and
	$\Re\,h_\IC(0)=\Im\,h_\IC(0)=0$.  We define $H=(H_1,\ldots,H_{2n})\in \IQ[Z]^{2n}$
	by letting $H_i$ to be the $(i+1)/2^\mathrm{th}$-component of $\Re\,h_\IC$ if $i$ is odd and the $i/2^\mathrm{th}$-component of $\Im\,h_\IC$ if $i$ is even where $i\in \{1,\ldots,2n\}$.
	We have
	$$\Re\, f(h(X))_\IC=\Re\, f_\IC(H(Z)), 
	\Im\, f(h(X))_\IC=\Im\, f_\IC(H(Z)).$$
	Hence $\Re\, f(h(X))_\IC, \Im\, f(h(X))_\IC\in\ma{IA}_{2n}$.
	
	\vs{0.2cm}
	$(\mathrm{IV}_{\IC})$:
	Let $a=(a_1,\ldots,a_n)\in \ma{D}(f)\cap \IQ^n$ and consider 
	$b:=(a_1,0,a_2,0,\ldots,a_n,0)\in \IQ^{2n}$. Then $b$ belongs to $\ma{D}(\Re\,f_\IC)$ and $\ma{D}(\Im\,f_\IC)$. We have
	$\Re\,(f_a)_\IC=(\Re\,f_\IC)_b$ and $\Im\,(f_a)_\IC=(\Im\,f_\IC)_b$.
	Hence $\Re\,(f_a)_\IC, \Im\,(f_a)_\IC\in\ma{IA}_{2n}$.
	
	\vs{0.2cm}
	$(\mathrm{V}_{\IC})$: 
	Let $F:=\int f\,dX_n$. We have $\ma{D}(F)=\ma{D}(f)$ and $\ma{B}(F)=\ma{B}(f)$. 
	For $z\in \ma{B}(f)$ we have
	$$F(z)=z_n\int_0^1f(z',tz_n)dt.$$
	Let 
	$$G:\ma{B}(F)\to \IC, z\mapsto \int_0^1 f(z',tz_n).$$
	Then $G\in \ma{O}_n$. It is enough to show that $\Re\,G_\IC, \Im\,G_\IC$ belong to $\ma{IA}_{2n}$.
	Let $f_1:=\Re\,f_\IC, f_2:=\Im\,f_\IC$ and $G_1:=\Re\,G_\IC, G_2:=\Im\,G_\IC$.
	Let $\Delta:=\ma{B}(f)\cap \ma{D}(f_1)\cap \ma{D}(f_2)$. Let $j\in \{1,2\}$. Then $\Delta\subset \ma{D}(G_j)$.
	For $z\in \Delta$ we have
	$$G_j(z)=\int_0^1 f_j(z',tx_n,ty_n)\,dt.$$
	We denote the integrand of the integral by $\hat{f}_j$.
	Let $T$ be a new variable.
	By the inductive assumption, Axiom (P1), Corollary 2.6 and Property (III) we have that $\hat{f}_j(Z,T)\in \ma{IA}_{2n+1}$.
	Let $\hat{G}_j:=\int \hat{f}_j\,dT$. By Property (V) we have $\hat{G}_j\in \ma{IA}_{2n+1}$.
	Let $c:=(0,\ldots,0,1)\in \IQ^{2n+1}$. Then $c\in \ma{D}(\hat{G}_j)$.
	We obtain
	$G_j(Z)=(\hat{G}_j)_c(Z,0)$. We are done by Properties (III) and (IV).
	\hfill$\blacksquare$
	
	\vs{0.5cm}
	{\bf 2.8 Corollary} (Integration)
	
	\vs{0.1cm}
	{\it Let $n\in \IN$ and $f\in\ma{IA}_n$. Let $R=(R_1,\ldots,R_n)\in \ma{R}(f)$ and let $R':=(R_1,\ldots,R_{n-1})$. We set $D:=D^n(R)$ and $D':=D^{n-1}(R')$.
		Let $a\in \ik_{>0}$ be such that $(0,a)\in D$. Let
		$$g:D'\to \IR, g(x')=\int_0^af(x',\zeta)\,d\zeta.$$
		Then the germ of $g$ at the origin belongs to $\ma{IA}_{n-1}$.}
	
	\vs{0.1cm}
	{\bf Proof:}
	
	\vs{0.1cm}
	We can use the proof of Proposition 1.23. Note that Properties (III), (IV) and (V) are sufficient.
	\hfill$\blacksquare$
	
	\vs{0.5cm}
	{\bf 2.9 Corollary}
	
	\vs{0.1cm}
	{\it Weierstra\ss $ $ preparation and division hold for $\ma{IA}$.}

	\vs{0.1cm}
	{\bf Proof:}
	
	\vs{0.1cm}
	By Proposition 2.4 we have that $\pi\in \ik$.
	By Corollary 2.8 and the Properties (I) - (V) we can use the proof of Theorem 1.24, choosing there $\delta\in \IQ_{>0}^n$.
	\hfill$\blacksquare$
	
	\vs{0.5cm}
	{\bf 2.10 Corollary} (Composition)
	
	\vs{0.1cm}
	{\it $\ma{IA}$ is closed under composition.}
	
	\vs{0.1cm}
	{\bf Proof:}
	
	\vs{0.1cm}
	The proof of Proposition 1.15(3) does only require the Weierstra\ss $ $ division theorem. Hence we are done by Corollary 2.9.
	\hfill$\blacksquare$
	
	\vs{0.5cm}
	Establishing the following result gives finally Theorem 2.2.
	
	\vs{0.5cm}
	{\bf 2.11 Proposition}
	
	\vs{0.1cm}
	{\it $\ma{IA}$ is a convergent system.}
	
	\vs{0.1cm}
	{\bf Proof:}
	
	\vs{0.1cm}
	Let $n\in \IN$ and $f\in \ma{IA}_n$. Let $a=(a_1,\ldots,a_n)\in \ma{D}(f)\cap \ik^n$. We have to show that $f_a\in \ma{IA}_n$.
	By Axiom (P3) and Property (III) we can assume that $a_i\neq 0$ for all $i\in \{1,\ldots,n\}$. Consider
	$h:=(a_1X_1,\ldots, a_nX_n)\in \ik[X]^n$. By Corollary 2.6 we have that $h\in \ma{IA}^n$. Since $h(0)=0$ we have by Corollary 2.10 that $g(X):=f(h(X))\in \ma{IA}_n$.
	We have that $b:=(1,\ldots,1)\in \ma{D}(g)\cap \IQ^n$. Hence $g_b\in \ma{IA}_n$ by Property (IV).
	As above we have that $\tilde{h}:=(X_1/a_1,\ldots,X_n/a_n)\in \ma{IA}_n$ with $\tilde{h}(0)=0$.
	Since $f_a(X)=g_b(\tilde{h}(X))$ we obtain again by Corollary 2.10 that $f_a\in \ma{IA}_n$.
	\hfill$\blacksquare$
	
	\vs{0.5cm}
	{\bf 2.12 Corollary}
	
	\vs{0.1cm}
	{\it The system $\ma{IA}$ of integrated algebraic series is the smallest convergent Analysis system.}
	
	\vs{0.1cm}
	{\bf Proof:}
	
	\vs{0.1cm}
	Let $\ma{S}$ be a convergent Analysis system.
	By Theorem 1.24 and Corollary 1.17 we obtain that $\ma{A}$ is contained in $\ma{S}$. Hence $\mathrm{CAS}(\ma{A})$ is contained in $\ma{S}$. We are done by Theorem 2.2.
	\hfill$\blacksquare$

	\subsection{The Field of Integrated Algebraic Numbers}
	
	{\bf 2.13 Definition} (Integrated Algebraic Numbers)
	
	\vs{0.1cm}
	We denote the coefficient field of $\ma{IA}$ by $\II\IA$ and call its elements {\bf integrated algebraic numbers}.
	
	\vs{0.5cm}
	We obtain Theorem A.
	
	\vs{0.5cm}
	{\bf 2.14 Corollary}
	
	\vs{0.1cm}
	{\it The field of integrated algebraic numbers is a countable real closed field.}
	
	\vs{0.1cm}
	{\bf Proof:}
	
	\vs{0.1cm}
	We have that $\ma{IA}$ is CAS by Theorem 2.2 and hence CWS by Theorem 1.24.
	By Theorem 1.16 we get that $\II\IA$ is real closed. 
	By Example 1.8 we have that the system $\ma{A}$ is countable. By Theorem 2.2 and Proposition 1.28 we have that $\ma{IA}$ and therefore $\II\IA=\ik(\ma{IA})$ is countable.
	\hfill$\blacksquare$

	\section{Convergent Systems and O-Minimality}
	
	\subsection{O-Minimality of Convergent Weierstra\ss $ $ Systems}
	
	We fix a convergent Weierstra\ss $ $ 
	system $(\ma{S}_n)_{n\in \IN_0}$ with coefficient field $\ik:=\ik(\ma{S})$.
	We set $I:=\{x\in \IR\mid -1\leq x\leq 1\}$ and $I_\ik:=I\cap \ik$.
	
	\vs{0.5cm}
	{\bf 3.1 Definition} (Restricted Function)
	
	\vs{0.1cm}
	Let $n\in \IN_0$. 
	A function $f:\ik^n\to \ik$ is called a {\bf restricted $\ma{S}$-function} if there is $p\in \ma{S}_n$ which converges on a neighbourhood of the unit cube $I^n$ 
	such that 
	$$f(x)=\left\{\begin{array}{ccc}
		p(x)&&x\in I_\ik^n,\\
		&\mbox{ if }&\\
		0&&x\notin I_\ik^n.
	\end{array}\right.$$
	By $\mathrm{R}\ma{S}_n$ we denote the set of restricted $\ma{S}$-functions of arity $n$. 
	
	\vs{0.5cm}
	Note that by Remark 1.9 the restricted functions are well-defined. The sets $\mathrm{R}\ma{S}_n$ are $\ik$-vector spaces which are closed under multiplication.
	
	\vs{0.5cm}
	Let $\mathrm{R}\ma{S}:=\bigcup_{n\in\IN_0}\mathrm{R}\ma{S}_n$ and let $\widetilde{\ma{S}}$ be the expansion of the coefficient field $\ik$ generated by $\mathrm{R}\ma{S}$. Note that in particular the graphs of unrestricted addition and multiplication on the field $\ik$ belong to $\widetilde{\ma{S}}$.  Let $\ma{L}$ be the language of ordered rings. Let $\ma{L}_{\ma{S}}$ denote the extension of $\ma{L}$ by symbols for every function from $\mathrm{R}\ma{S}$. By $T_\ma{S}$ we denote the $\ma{L}_\ma{S}$-theory of $\widetilde{\ma{S}}$. Let $\ma{L}_\ma{S}^{-1}$ be the extension of $\ma{L}_\ma{S}$ by a symbol for the reciprocal $\ik_{>0}\to \ik_{>0}, x\mapsto 1/x$.  Let $\ma{L}_\ma{S}^\IQ$ be the extension of $\ma{L}_\ma{S}$ by symbols for every power function $\ik_{>0}\to \ik_{>0}, x\mapsto x^q,$ where $q\in \IQ$. By Theorem 1.16 the power functions take values in $\ik$ and are therefore well-defined. Moreover, we view $\widetilde{\ma{S}}$ in the natural way as an $\ma{L}_\ma{S}^\IQ$-structure (extending the power functions by $0$ to $\ik_{\leq 0}$). By 'definable' we mean definable in the $\ma{L}_\ma{S}$-structure $\widetilde{\ma{S}}$. 
	Note that we do not have to care about parameters sind the universe is contained in the language via $\mathrm{R}\ma{S}_0$.
	
	\vs{0.5cm}
	{\bf 3.2 Theorem}
	
	\vs{0.1cm}
	{\it The structure $\widetilde{\ma{S}}$ has quantifier elimination in the language $\ma{L}_\ma{S}^{-1}$.}
	
	\vs{0.1cm}
	{\bf Proof:}
	
	\vs{0.1cm}
	We want to establish [10, Proposition 2.9] for the field $\ik$ as there, adopting the proof of [6, Theorem 4.6]. The condition of [10, Proposition 2.9] that the field $\ik$ is real closed is fulfilled by Theorem 1.16. We have all the arguments involving Weierstra\ss $ $ preparation in hand. As discussed in [10, (2.10)] the key step to adopt the proof of [6, Theorem 4.6] is to establish a covering argument. This requires extra work since $\ik$ is not locally compact for $\ik\neq\IR$. The covering argument occurs in the proofs of [6, Basic Lemma 4.10 \& Local Basic Lemma 4.13]. Inspecting these proofs we are done by the following lemma and Axiom (C).
	\hfill$\blacksquare$
	
	\vs{0.5cm}
	We need some notations.
	Fix $n\in \IN$ and let $f\in \ma{O}_n$ be regular in $X_n$.
	Let $\ma{D}:=\ma{D}(f)$ and $\ma{D}'$ be the projection of $\ma{D}$ onto the first $n-1$ coordinates. 
	We set
	$$\Sigma:=\{a\in \ma{D}'\mid f_{(a,0)}\mbox{ is regular in }X_n\}.$$
	For $a\in \Sigma$ we let $d(a)\in \IN_0$ be the order of $f_{(a,0)}$ in $X_n$.
	For $a\in \Sigma$ let $P_a\in \ma{O}_{n-1}[X_n]$ be the unique Weierstra\ss $ $ polynomial in $X_n$ of degree $d(a)$ and $u_a\in \ma{O}_n$ be the unique power series with $u_a(0)\neq 0$ such that $f_{(a,0)}=P_au_a$.
	Moreover, we set for $a\in \Sigma$
	$$\ma{C}_a:=\ma{D}(P_a)\cap\{x\in \ma{D}(u_a)\mid \mathrm{sign}(u_a(x))=\mathrm{sign}(u_a(0))\}$$
	and let $\ma{C}_a'$ be the projection of $\ma{C}_a$ to the first $n-1$ coordinates. Then $(\ma{C}_a'+a)_{a\in \Sigma}$ is obviously a covering of $\Sigma$ by open sets.

	Finally, for $a\in \IR^{n-1}$ we set $\ma{Z}(a):=\{i\in \{1,\ldots,n-1\}\mid a_i= 0\}$.
	
	\vs{0.5cm}
	{\bf 3.3 Lemma}
	
	\vs{0.1cm}
	{\it Assume  that $\Sigma$ is open and $d(a)=d(b)$ for every $a,b\in \Sigma$ with $\ma{Z}(a)=\ma{Z}(b)$. Then $(\ma{C}_a'+a)_{a\in \Sigma\cap \IQ^{n-1}}$ is a covering of $\Sigma$.}
	
	\vs{0.1cm}
	{\bf Proof:}
	
	\vs{0.1cm}
	Let $a\in \Sigma$. Since $\Sigma$ is open we find some $R\in \IQ_{>0}^{n-1}$ such that $D^{n-1}(2R)+a\subset \Sigma$.
	Let $\Delta:=(D^{n-1}(R)\cap \ma{C}_a')+a$.
	For $b\in \Delta$ we have that 
	$$f=P_a(X'-a,X_n)u_a(X'-a,X_n)=P_b(X'-b,X_n)u_b(X'-b,X_n)$$
	on the non-empty open subset $(\ma{C}_a+(a,0))\cap (\ma{C}_b+(b,0))$. 
	If $d(a)=d(b)$ we obtain that $P_b=(P_a)_{(b-a,0)}, u_b=(u_a)_{(b-a,0)}$ and therefore $\Delta\subset \ma{C}_b'+b$. 
	Choosing $b\in \Delta\cap \IQ^{n-1}$ with $\ma{Z}(b)=\ma{Z}(a)$ we are done be the assumption.
	\hfill$\blacksquare$
	
	\vs{0.5cm}
	{\bf 3.4 Corollary}
	
	\vs{0.1cm}
	{\it The structure $\widetilde{\ma{S}}$ is o-minimal.}
	
	\vs{0.1cm}
	{\bf Proof:}
	
	\vs{0.1cm}
	This follows from Theorem 3.2 and the following claim.
	
	\vs{0.2cm}
	{\bf Claim:} Let $f\in \ma{S}_1$ be with $f\neq 0$ and let $a\in \ma{D}(f)$ be with $f(a)=0$. Then $a\in \ik$.
	
	\vs{0.1cm}
	{\bf Proof of the claim:}
	One can adjust the proof of Theorem 1.16.
	\hfill$\blacksquare_\mathrm{Claim}$
	
	\hfill$\blacksquare$

	\vs{0.5cm}
	Let $\IR\widetilde{\ma{S}}$ be the canonical $\ma{L}_\ma{S}$ structure on the reals (evaluating $p\in \ma{S}_n$ converging on an open neighbourhood of $I^n$ there and not just $I_\ik^n$). 
	
	\vs{0.5cm}
	{\bf 3.5 Proposition}
	
	\vs{0.1cm}
	{\it $\IR\widetilde{\ma{S}}$ is the only model of $T_\ma{S}$ on the reals extending $\widetilde{\ma{S}}$.}
	
	\vs{0.1cm}
	{\bf Proof:}
	
	\vs{0.1cm}
	Uniqueness follows easily from the fact that $\ik$ is dense in $\IR$ and continuity arguments.
	For existence we consider $\IR\widetilde{\ma{S}}$ canonically as an $\ma{L}_\ma{S}^{-1}$-structure. 
	The proof of Theorem 3.2 can be taken to show that also $\IR\widetilde{\ma{S}}$ has quantifier elimination in $\ma{L}_\ma{S}^{-1}$.
	Hence $\IR\widetilde{\ma{S}}$ is elementarily equivalent to $\widetilde{\ma{S}}$ as $\ma{L}_\ma{S}^{-1}$ structures and therefore as $\ma{L}_\ma{S}$-structures.
	\hfill$\blacksquare$
	
	\vs{0.5cm}
	{\bf 3.6 Theorem}
	
	\vs{0.1cm}
	{\it The structure $\widetilde{\ma{S}}$ has universal axiomatization in the language $\ma{L}_\ma{S}^\IQ$.}
	
	\vs{0.1cm}
	{\bf Proof:}
	
	\vs{0.1cm}
	One can adjust the reasoning of [10, Section 2] to establish that $\IR\widetilde{\ma{S}}$ has universal axiomatization in the language $\ma{L}_\ma{S}^\IQ$.
	The universe in our setting is also the field of reals, but we are only working with power series from $\ma{S}$. This can be done since the properties of a Weierstra\ss $ $ system (see Proposition 1.15) are sufficient. Note that the equivalent of [10, Corollary 2.11] holds by Theorem 3.2 and Proposition 3.5 (see also Rambaud [30]).
	Since $\widetilde{\ma{S}}$ is elementarily equivalent to $\IR\widetilde{\ma{S}}$ we obtain that $\widetilde{\ma{S}}$ has universal axiomatization in the language $\ma{L}_\ma{S}^\IQ$.
	\hfill$\blacksquare$
	
	\vs{0.5cm}
	{\bf 3.7 Corollary}
	
	\vs{0.1cm}
	{\it Functions definable in $\widetilde{\ma{S}}$ are piecewise given by $\ma{L}_\ma{S}^\IQ$-terms.}
	
	\vs{0.1cm}
	{\bf Proof:}
	
	\vs{0.1cm}
	This follows by quantifier elimination Theorem 3.2 and universal axiomatization Theorem 3.6 (compare the proof of [10, Corollary 2.11]).
	\hfill$\Box$
	
	\vs{0.5cm}
	{\bf 3.8 Remark}
	
	\vs{0.1cm}
	That functions definable in $\widetilde{\ma{S}}$ are piecewise given by $\ma{L}_\ma{S}^\IQ$-terms can be expressed in the following way:

	For $n\in\IN_0$ let $\Omega_n$ be the union of
	$\ik[X]$ and $\mathrm{R}\ma{S}_n$ and, if $n=1$, of all power functions $\ik_{>0}\to \ik_{>0}, x\mapsto x^q,$ with rational exponent $q$, extended by $0$ to $\ik_{\leq 0}$. Note that these functions are well-defined by the fact that $\ik$ is real closed.
	Finally we set $\Omega:=\bigcup_{n,l\in\IN}\Omega_n^l$.
	
	Let $f:\ik^n\to \ik$ be definable in $\widetilde{\ma{S}}$. Then there is a partition of $\ik^n$ into finitely many sets definable in $\widetilde{\ma{S}}$ such that $f$ is a finite composition of elements from $\Omega$ on each of these sets.
	
	\vs{0.5cm}
	{\bf 3.9 Comment}
	
	\vs{0.1cm}
	Rambaud [30] has established the above results for convergent Weierstra\ss $ $ systems on the reals. There one can see in detail how the adjustments of the reasoning in [6] and [10] to this setting is done.
	Here we have the additional difficulty that we are not on the reals. This is handled by our Lemma 3.3.

	\vs{0.5cm}
	We establish Lion-Rolin preparation [26] for $\widetilde{\ma{S}}$ (see also [28] and [19]).
	
	\vs{0.5cm}
	{\bf 3.10 Definition} (Lion-Rolin Preparation)
	
	\vs{0.1cm}
	We say that $\widetilde{\ma{S}}$ has {\bf Lion-Rolin preparation} if the following holds.
	Let $A\subset \ik^n$  be definable in $\widetilde{\ma{S}}$  and let $f:A\to \ik, (x',x_n)\to f(x',x_n),$ be a definable function. Then there is a cell decomposition $\ma{C}$ of $A$ such that the following holds. Let $C\in\ma{C}$ and let $B$ denote the base of $C$; i.e. $B$ is the projection of $C$ onto the first $n-1$ coordinates.
	Assume that $C$ is fat with respect to the last variable $x_n$; i.e. $C_{x'}$ is a nonempty open interval (and not just a point) for every $x'\in B$.
	Then the function $f|_C$ can be written as
	$$f|_C(x',x_n)=a(x')\big\vert x_n-\theta(x')\big\vert ^ru\big(x',x_n-\theta(x')\big)$$
	where
	$r\in \IQ$, the functions $a, \theta:B\to \ik$ are definable in $\widetilde{\ma{S}}$, $x_n\neq \theta(x')$ on $C$, and $u(x',x_n)$ is a so-called special unit on
	$$C^\theta:=\big\{\big(x',x_n-\theta(x')\big)\mid (x',x_n)\in C\big\};$$
	i.e. $u$ is of the form
	$$u(x',x_{n-1})=v\Big(b_1(x'),\ldots,b_M(x'),b_{M+1}(x')|x_n|^{1/q},b_{M+2}(x')|x_n|^{-1/q}\Big)$$
	where $q\in\IN$,
	\begin{eqnarray*}
		\varphi:B\times \ik\setminus\{0\} &\to&\ik^{M+2},\\
		(x',x_n)&\mapsto&\Big(b_1(x'),\ldots,b_M(x'),b_{M+1}(x')|x_n|^{1/q},b_{M+2}(x')|x_n|^{-1/q}\Big),
	\end{eqnarray*}
	is a definable function with $\varphi(C^\theta)\subset [-1,1]^{M+2}$ and $v\in \mathrm{R}\ma{S}_{M+2}$ with $v([-1,1]^{M+2})\subset [1/c,c]$ for some $c\in \ik$ with $c>1$.
	
	\vs{0.5cm}
	{\bf 3.11 Theorem}
	
	\vs{0.1cm}
	{\it The structure $\widetilde{\ma{S}}$ has Lion-Rolin preparation.}
	
	\vs{0.1cm}
	{\bf Proof:}
	
	\vs{0.1cm}
	Since we already know that $\widetilde{\ma{S}}$ is o-minimal we can use in the arguments below cell decomposition.
	By Corollary 3.7 definable functions are piecewise given by terms.
	By passing to a suitable cell decomposition we may therefore assume that 
	the definable function is given by a term.
	We do induction on the complexity of the term.
	In the base case we deal with the constants, coordinate functions, addition, multiplication, power functions and restricted $\ma{S}$-functions.
	Constants, coordinate functions and power functions are already prepared. Addition and multiplication can be handled easily. For restricted $\ma{S}$-functions we can adopt the proof of [26, Proposition 2]. There Weierstra\ss $ $ preparation is used and the finiteness result [6, Lemma 4.12]. Note that the algebraic arguments involved there are valid in our setting by Remark 1.14.
	For the inductive step we have to deal by the base case with compositions of functions which can be prepared. We can use [28, Lemma 4.4 \& Lemma 5.3].
	\hfill$\blacksquare$
	
	\vs{0.5cm}
	{\bf 3.12 Comment}
	
	\begin{itemize}
		\item[(1)]
		The orginal Lion-Rolin preparation in the real case in [26] includes that the involved cells are real analytic and that the functions depending only on $x'$ (i.e. the corresponding $a(x'),
		\Theta(x'), b_1(x'),\ldots,b_{M+2}(x')$) are real analytic. We could also obtain this by some extra technical effort by defining the canonical extension of functions to the reals (compare with the beginning of Section 4) and demanding that these are real analytic. But we refrain from doing so since it is not necessary for the desired results on integration. 
		\item[(2)]
		In [28] Lion-Rolin preparation with additional information on the cells is established for convergent systems on the reals with Weierstra\ss $ $ preparation being closed under composition (but not necessarily having Weierstra\ss $ $ division).
		As stated there the full Lion-Rolin preparation implies quantifier elimination, o-minimality and the description of definable functions by terms. We have done it the other way round.
		As observed in [28, Footnote on p. 4396] one obtains from our starting point immediately a weaker version of preparation as in van den Dries and Speissegger [11]. 
		To establish our Lion-Rolin preparation we only had to invest some work.
		For establishing full Lion-Rolin preparation as in [28] one has also to care inductively about cell decomposition which requires heavy analysis.
	\end{itemize}
	
	\subsection{Convergent Analysis Systems and Exponentiation}
	
	We fix a convergent Analysis system $\ma{S}$ with coefficient field $\ik:=\ik(\ma{S})$.
	
	\vs{0.5cm}
	{\bf 3.13 Proposition}
	
	\vs{0.1cm}
	{\it The field $\ik$ is closed under logarithm and exponentiation.}
	
	\vs{0.1cm}
	{\bf Proof:}
	
	\vs{0.1cm}
	We start with the logarithm. Let $a\in \ik_{>0}$. We have to show that $\log(a)\in \ik$.
	By the properties of the logarithm we can assume that $a\leq 1$. We have that $f:=1/(1+X)\in \ma{S}_1$ by Axioms (P2) and (A1).
	Hence by Axiom (A5) the logarithmic series $L(X)=\int f\,dX$ belongs to $\ma{S}_1$.
	We obtain by Remark 1.9 that
	$\log(a)=L(a-1)\in\ik$.
	
	\vs{0.2cm}
	For exponentiation we have to show that $\exp(a)\in \ik$ for $a\in \ik$.
	Let $E(X)=\sum_{p=0}^\infty X^p/p!$ be the exponential series.
	Since $E-1$ is the inverse of the logarithmic series we obtain by Theorem 1.24 and Proposition 1.15(4) that $E\in \ma{S}_1$.
	By Remark 1.9 we get that $\exp(a)=E(a)\in \ik$.
	\hfill$\blacksquare$
	
	\vs{0.5cm}
	Let $\ma{S}(\exp)$ be the expansion of $\ma{S}$ by the exponential function $\exp:\ik\to \ik_{>0}$ on $\ik$.
	Let $\ma{L}_\ma{S}(\exp,\log)$ be the extension of the language $\ma{L}_\ma{S}(\exp)$ by a symbol for the exponential function on $\ik$ and by a symbol for the logarithmic function $\log:\ik_{>0}\to \ik$ on $\ik$ (extended by $0$ to the negative axis).
	We view $\ma{S}(\exp)$ in the natural way as an $\ma{L}_\ma{S}(\exp,\log)$-structure.
	
	\vs{0.5cm}
	{\bf 3.14 Theorem}
	
	\vs{0.1cm}
	{\it The structure $\widetilde{\ma{S}}(\exp)$ has quantifier elimination in the language $\ma{L}_\ma{S}(\exp,\log)$ and has universal axiomatization there.}
	
	\newpage
	{\bf Proof:}
	
	\vs{0.1cm}
	We can take the reasoning in [10, Sections 3 \& 4], using the result of the Section 3.1 above. For the necessary adjustments see [10, Section (3.9)].
	\hfill$\blacksquare$
	
	\vs{0.5cm}
	{\bf 3.15 Corollary}
	
	\vs{0.1cm}
	{\it Functions definable in $\widetilde{\ma{S}}(\exp)$ are piecewise given by $\ma{L}_\ma{S}(\exp,log)$-terms.}
	
	\vs{0.5cm}
	Let $\IR\widetilde{\ma{S}}(\exp)$ be the canonical $\ma{L}_\ma{S}(\exp,\log)$-structure on the reals.
	
	\vs{0.5cm}
	{\bf 3.16 Corollary}
	
	\vs{0.1cm}
	{\it The structure $\widetilde{\ma{S}}(\exp)$ is o-minimal.}
	
	\vs{0.1cm}
	{\bf Proof:}
	
	\vs{0.1cm}
	The structure $\IR\widetilde{\ma{S}}(\exp)$ has the same universal axiomatization as $\widetilde{\ma{S}}(\exp)$ in the language $\ma{L}_\ma{S}(\exp,\log)$. One has to do the arguments of the proof of Theorem 3.14 twice (see also [30]). 
	Hence the $\ma{L}_\ma{S}(\exp,\log)$-structures $\widetilde{\ma{S}}(\exp)$ and $\IR\widetilde{\ma{S}}(\exp)$ are elementarily equivalent.
	The latter is a reduct of the o-minimal structure $\IR_{\an,\exp}$ and therefore o-minimal. Since o-minimality is a first order property we are done.
	\hfill$\blacksquare$
	
	\vs{0.5cm}
	{\bf 3.17 Comment}
	
	\vs{0.1cm}
	Rambaud [30] has established the above results for convergent Weierstra\ss $ $ systems on the reals. There one can see in detail how the adjustments of the reasoning in [10] to this setting is done.

	\section{Convergent Analysis Systems and Integration}
	
	Let $\mathbb{X}$ be a dense subset of the reals and let $\ma{R}$ be an o-minimal structure on $\mathbb{X}$. By cell decomposition and induction, starting with intervals and points, one defines for every set $A\subset \mathbb{X}^n$ or function $f:A\to \mathbb{X}$ definable in $\ma{R}$ a canonical extension $A_\IR\subset \IR^n$ respectively $f_\IR:A_\IR\to \IR$ 
	(with $A_\IR\cap \mathbb{X}^n=A$ and $f_\IR|_A=f$.)
	The extended sets form an o-minimal structure on $\IR$ denoted by $\ma{R}_\IR$.
	
	\vs{0.5cm}
	{\bf 4.1 Example}
	
	\begin{itemize}
		\item[(1)]
		Given a convergent Weierstra\ss $ $ system $(\ma{S}_n)_{n\in \IN_0}$ with coefficient field $\ik$ we have $\widetilde{\ma{S}}_\IR=\IR\widetilde{\ma{S}}$.
		Given a formula $\varphi(X)$ in the language $\ma{L}_\ma{S}$ we have $\varphi(\ik)_\IR=\varphi(\IR)$.
		\item[(2)]
		Given a convergent Analysis system $(\ma{S}_n)_{n\in \IN_0}$ with coefficient field $\ik$ we have $\widetilde{\ma{S}}(\exp)_\IR=\IR\widetilde{\ma{S}}(\exp)$.
		Given a formula $\varphi(X)$ in the language $\ma{L}_\ma{S}(\exp,\log)$ we have $\varphi(\ik)_\IR=\varphi(\IR)$.
	\end{itemize}
	
	\vs{0.2cm}
	We fix an o-minimal structure $\ma{R}$ on a dense subset $\mathbb{X}$ of $\IR$. 
	
	\vs{0.5cm}
	{\bf 4.2 Definition}
	
	\vs{0.1cm}
	Let $n\in \IN$ and $f:\mathbb{X}^n\to \mathbb{X}$ be definable in $\ma{R}$. We say that $f$ is {\bf integrable} if $f_\IR:\IR^n\to \IR$ is (Lebesgue-)integrable. In that case we set
	$$\int_{\mathbb{X}^n} f(x)\,dx:=\int_{\IR^n}f_\IR(x)\, dx.$$
	
	\vs{0.5cm}
	Note that by cell decomposition the function $f_\IR$ is measurable (see [18, Proposition 1.1]). The function $f$ is integrable if and only if 
	$\int_{\IR^n} |f_\IR(x)|\,dx<\infty$ (i.e. if and only if $f$ is absolutely integrable). It is no loss of generality that we assume that $f$ is defined on $\mathbb{X}^n$. If $f:A\to \mathbb{X}$ is definable in $\ma{R}$ where $A\subset  \mathbb{X}^n$ (is definable in $\ma{R}$) we may always extend it to $\mathbb{X}^n$ by zero.
	Note that the integral could be also defined internally (see [17]).
	
	\vs{0.5cm}
	{\bf 4.3 Definition}
	
	\vs{0.1cm}
	Let $f:\mathbb{X}^m\times \mathbb{X}^n\to \mathbb{X}, (x,y)\to f(x,y),$ be definable in $\ma{R}$. We set
	$$\mathrm{Fin}(f):=\big\{x\in \mathbb{X}^m\;\big\vert\; f_x\mbox{ is integrable}\big\}$$
	and 
	$$\mathrm{Int}(f):\mathrm{Fin}(f)\to \IR, x\mapsto \int_{\mathbb{X}^n} f(x,y)\,dy.$$
	
	\vs{0.5cm}
	From now on we assume that $\ma{R}$ expands a subfield $\ik$ of the reals.
	
	\vs{0.5cm}
	{\bf 4.4 Definition} (Volume Ring)
	
	\vs{0.1cm}
	We call
	$$\mathrm{Vol}(\ma{R}):=\Big\{\int_{\ik^n}f(x)\,dx\;\big\vert\; n\in \IN_0, f:\ik^n\to \ik\mbox{ definable in } \ma{R} \mbox{ and integrable}\Big\}$$
	the {\bf volume ring} of $\ma{R}$.
	
	\vs{0.5cm}
	{\bf 4.5 Remark}
	
	\vs{0.1cm}
	$\mathrm{Vol}(\ma{R})$ is a $\ik$-subalgebra of $\IR$.
	
	\newpage
	{\bf 4.6 Definition} (Volume Closed)
	
	\vs{0.1cm}
	We call $\ma{R}$ {\bf volume closed} if $\mathrm{Vol}(\ma{R})=\ik$.
	
	\vs{0.5cm}
	From now on we fix a convergent Analysis system $\ma{S}=(\ma{S}_n)_{n\in \IN_0}$ with coefficient field $\ik$.
	
	\vs{0.5cm}
	{\bf 4.7 Proposition}
	
	\vs{0.1cm}
	{\it Let $f:\ik^m\times \ik^n\to \ik, (x,y)\to f(x,y)$ be definable in $\widetilde{\ma{S}}$. Then the set 
		$\mathrm{Fin}(f)$
		is definable in $\widetilde{\ma{S}}$.}
	
	\vs{0.1cm}
	{\bf Proof:}
	
	\vs{0.1cm}
	Since $\IR\widetilde{\ma{S}}$ is a reduct (in the sense of definability) of the polynomially bounded o-minimal structure $\IR_\an$ (see [10] for this structure)  it is polynomially bounded. Hence $\widetilde{\ma{S}}$ is polynomially bounded.
	We can adopt the proof of [18, Theorem 2.2b)].
	\hfill$\blacksquare$

	\vs{0.5cm}
	{\bf 4.8 Theorem}
	
	\vs{0.1cm}
	{\it Let $f:\ik^m\times\ik^n\to \ik$ be definable in $\widetilde{\ma{S}}$. Then there are finite numbers  $K,L\in\IN_0$ and for $k\in \{1,\ldots,K\}$ and $l\in \{1,\ldots,L\}$ there are functions $\varphi_k:\mathrm{Fin}(f)\to \ik, \psi_{kl}:\mathrm{Fin}(f)\to \ik_{>0}$ definable in $\widetilde{\ma{S}}$  such that $$\mathrm{Int}(f)=\sum_{k=1}^L\varphi_k\Big(\prod_{l=1}^L\log(\psi_{kl})\Big).$$}
	
	\vs{0.1cm}
	{\bf Proof:}
	
	\vs{0.1cm}
	The system $\ma{S}$ is a convergent Weierstra\ss $ $ system by Theorem 1.24. It is closed under reciprocals by Axiom (A1) and antidifferentiation by Axiom (A5). By Proposition 1.15(2) it is closed under composition and by Axiom (A4) it is closed under complexification. Therefore as observed in [28, Proposition 9.3] the Lion-Rolin splitting (see [27, I.6]) holds:
	
	Let $f\in\ma{S}_{n+1}$. Then there are $f_+,f_-\in\ma{S}_{n+1}$ such that for all sufficiently small $x,y/z$ and $z\neq 0$
	$$f(x,y/z,z)=f_+(x,y,z)+(y/z)f_-(x,y,y/z).$$
	
	With this in hand the techniques (where formal antiderivatices are taken) of [26] and [5] can be literally translated to our situation and we are done.
	\hfill$\blacksquare$
	
	\vs{0.5cm}
	{\bf 4.9 Corollary}
	
	\vs{0.1cm}
	{\it The structure $\widetilde{\ma{S}}$ is volume closed.}
	
	\vs{0.1cm}
	{\bf Proof:}
	
	\vs{0.1cm}
	This follows by Theorem 4.8 (in the case $m=0$) and Proposition 3.13.
	\hfill$\blacksquare$
	
	\section{Integrated Algebraic Numbers and Periods}
	
	We consider the following three o-minimal structures: $\widetilde{\ma{A}}, \widetilde{\ma{IA}}$ and $\widetilde{\ma{IA}}(\exp)$. 
	Note that $\widetilde{\ma{A}}$ is the pure field structure $\IA$, the definable sets are the sets semialgebraic over the rationals.
	
	\vs{0.3cm}
	Recall that a {\bf period} is the integral of an integrable function which is semialgebraic over the rationals.
	By $\mathbb{P}$ we denote the ring of periods.
	
	\vs{0.5cm}
	{\bf 5.1 Proposition}
	
	{\it \begin{itemize}
			\item[(1)]
			The field $\IA$ is not volume closed. We have $\mathrm{Vol}(\IA)=\mathbb{P}$.
			\item[(2)]
			The structure $\widetilde{\ma{IA}}$ is volume closed.
	\end{itemize}}
	{\bf Proof:}
	
	\vs{0.1cm}
	(1): That $\mathrm{Vol}(\IA)=\mathbb{P}$ follows from the definition of periods. The first statement follows from the fact that 
	the circle number $\pi$ is a transcendental period.
	
	\vs{0.2cm}
	(2): This follows from Theorem 2.2 and Corollary 4.9.
	\hfill$\blacksquare$
	
	\vs{0.5cm}
	We obtain Theorem B.
	
	\vs{0.5cm}
	{\bf 5.2 Corollary}
	
	\vs{0.1cm}
	{\it A period is an integrated algebraic number.}
	
	\vs{0.1cm}
	{\bf Proof:}
	
	\vs{0.1cm}
	By Proposition 5.1(1) we know that $\mathbb{P}=\mathrm{Vol}(\widetilde{\ma{A}})$. Since $\ma{IA}$ contains $\ma{A}$ we have $\mathrm{Vol}(\widetilde{\ma{A}})\subset \mathrm{Vol}(\widetilde{\ma{IA}})$. By Proposition 5.1(2) we know that $\mathrm{Vol}(\widetilde{\ma{IA}})=\ik(\ma{IA})$ which is $\II\IA$.
	Combining these facts we obtain that $\mathbb{P}\subset \II\IA$. 
	\hfill$\blacksquare$
	
	\vs{0.5cm}
	In number theory often the {\bf extended period ring} $\mathbb{P}[1/\pi]$ is considered. Since $\pi\in \mathbb{P}$ and $\II\IA$ is a field we obtain that $\mathbb{P}[1/\pi]\subset \II\IA$.
	
	\vs{0.5cm}
	In our setting we can also deal with families of periods.
	
	\vs{0.5cm}
	{\bf 5.3 Definition} (Parameterized Periods)
	
	\vs{0.1cm}
	A {\bf parameterized period (of arity $m\in \IN_0$)} is a function $P: \IA^m\to \IR$ such that the following holds.
	There is a partition $\ma{C}$ of $\IA^m$ into finitely many sets which are semialgebraic over $\IA$ (or equivalently over $\IQ$) and for every $C\in \ma{C}$ there is $n_C\in \IN_0$ and a function $f_C:\IA^m\times\IA^{n_C}\to \IA$ semialgebraic over $\IA$ (or equivalently over $\IQ$) such that $C\subset \mathrm{Fin}(f_C)$ and $P(x)=\mathrm{Int}(f_C)(x)$ for all $x\in C$.
	
	\newpage
	{\bf 5.4 Remark}
	
	\begin{itemize}
		\item[(a)] A period is a parameterized period of arity $0$. 
		\item[(b)]
		Let $P: \IA^m\to \IR$ be a parameterized period and let $a\in \IA^m$. Then $P(a)$ is a period.
	\end{itemize} 
	
	\vs{0.2cm}
	{\bf 5.5 Theorem}
	
	\vs{0.1cm}
	{\it Let $P:\IA^m\to \IA$ be a parameterized period. 
		Then there are finite numbers $K,L\in\IN_0$ and for $k\in \{1,\ldots,K\}$ and $l\in \{1,\ldots,L\}$ there are functions $\varphi_k:\II\IA^m\to \II\IA, \psi_{kl}:\II\IA^m\to \II\IA_{>0}$ definable in $\widetilde{\ma{IA}}$ such that 
		$$P(x)=\sum_{k=1}^K\varphi_k(x)\Big(\prod_{l=1}^L\log\big(\psi_{kl}(x)\big)\Big)$$
		for all $x\in \IA^m$.}
	
	\vs{0.1cm}
	{\bf Proof:}
	
	\vs{0.1cm}
	For $C\in \ma{C}$ let $C_{\II\IA}$ and $f_{C,\II\IA}$ be the canonical extension of $C$ and $f_C$ to the real closed field $\II\IA$, respectively (see [2, Chapter 5]).
	Note that $\mathrm{Fin}(f_C)=\mathrm{Fin}(f_{C,\II\IA})\cap \IA^m$. We apply Theorem 4.8 and are done.
	\hfill$\blacksquare$
	
	\vs{0.5cm}
	We extend the notion of integrated algebraic numbers by the structure $\widetilde{\ma{S}}(\exp)$.
	
	\vs{0.5cm}
	{\bf 5.6 Definition} (Exponential integrated algebraic number)
	
	\vs{0.1cm}
	We set
	$\II\IA_\exp:=\mathrm{Vol}(\widetilde{\ma{IA}}(\exp))$. An element of $\II\IA_\exp$ is called an {\bf exponential integrated algebraic number}.
	
	\vs{0.5cm}
	We establish Theorem C.
	
	\vs{0.5cm}
	{\bf 5.7 Proposition}
	
	\vs{0.1cm}
	{\it $\II\IA_\exp$ is a countable $\II\IA$-subalgebra of $\IR$.}
	
	\vs{0.1cm}
	{\bf Proof:}
	
	\vs{0.1cm}
	By Theorem 2.2 and Proposition 1.28 we have that $\ma{IA}$ is countable. This implies that the language $\ma{L}_{\ma{IA}}$ and therefore the language $\ma{L}_{\ma{IA}}(\exp)$ is countable. Since
	the universe $\II\IA$ of the structure $\widetilde{\ma{IA}}(\exp)$ is countable by Corollary 2.14 we obtain that there are only countably many $\widetilde{\ma{IA}}(\exp)$-definable sets and functions. This gives the proposition.
	\hfill$\blacksquare$
	
	\vs{0.5cm}
	An {\bf exponential period} in the sense of Kontsevich and Zagier [25] is given by $\int_{\IR^n} \exp(f(x))g(x)\,dx$ with absolutely integrable integrand where $f=f_1+if_2$ and $g=g_1+ig_2$ (with functions $f_1,f_2,g_1,g_2:\IR^n\to \IR$) are semialgebraic over $\IQ$ (see also [14, 29]).
	In Commelin, Habegger and Huber [4] the additional condition on the function $f$ having bounded imaginary part $f_2$ is imposed.
	Their definition is sufficient for arithmetic geometry and we adopt it here.
	By $\mathbb{P}_\exp$ we denote the ring of exponential periods. 
	
	\vs{0.5cm}
	{\bf 5.8 Proposition}
	
	\vs{0.1cm}
	{\it The field of integrated algebraic numbers is closed under sine and cosine.}
	
	\vs{0.1cm}
	{\bf Proof:}
	
	\vs{0.1cm}
	We show this for the sine. For the cosine the arguments are the same.
	The arcsine series is an integrated algebraic series (see the introduction).
	By Theorem 2.2, Theorem 1.24 and Proposition 1.15(4) we get that the sine series is an integrated algebraic power series. It has radius of convergence $\infty$. We are done by Remark 1.9.
	\hfill$\blacksquare$
	
	\vs{0.5cm}
	{\bf 5.9 Proposition}
	
	\vs{0.1cm}
	{\it The restricted sine function
		$$\II\IA\to \II\IA, x\mapsto \left\{
		\begin{array}{ccc}
			\sin(x),&&x\in [-\pi,\pi]\cap\II\IA,\\
			&\mbox{ if }&\\
			0,&&x\notin [-\pi,\pi]\cap\II\IA,\\
		\end{array} \right.$$
		and the restricted cosine function
		$$\II\IA\to \II\IA, x\mapsto \left\{
		\begin{array}{ccc}
			\sin(x),&&x\in [-\pi,\pi]\cap\II\IA,\\
			&\mbox{ if }&\\
			0,&&x\notin [-\pi,\pi]\cap\II\IA,\\
		\end{array} \right.$$
		are well-defined and definable in $\widetilde{\ma{IA}}$.}
	
	\vs{0.1cm}
	{\bf Proof:}
	
	\vs{0.1cm}
	We have $\pi\in \II\IA$, for example by Corollary 5.2. We get by Proposition 5.8 that the functions are well-defined.
	By the proof of Proposition 5.8 and the definition of $\widetilde{\ma{IA}}$ the statement holds if we replace $[-\pi,\pi]$ by $I$. Using the angle addition theorems we obtain the result.
	\hfill$\blacksquare$

	\vs{0.5cm}
	{\bf 5.10 Proposition}
	
	\vs{0.1cm}
	{\it An exponential period is an exponential integrated algebraic number.}
	
	\vs{0.1cm}
	{\bf Proof:}
	
	\vs{0.1cm}
	This is a consequence of the chosen definition of an exponential period and the following observation. 
	Let $R\in \II\IA_{>0}$ and let $S:=\{z\in \IC\mid -R\leq \Im(z)\leq R\}$.
	The function  
	$$\II\IA(i)\to \II\IA(i), z\mapsto \left\{
	\begin{array}{ccc}
		\exp(z),&&z\in S\cap\II\IA(i),\\
		&\mbox{ if }&\\
		0,&&z\notin S\cap\II\IA(i),\\
	\end{array} \right.$$
	restricting the exponential function, is well-defined and definable in $\widetilde{\ma{IA}}(\exp)$.
	To show this note that
	$\exp(z)=\exp(\Re\,z)(\cos(\Im\,z)+i\sin(\Im\,z))$ for $z\in \IC$. 
	By Proposition 3.13 and Proposition 5.8 the function is well-defined. 
	Replacing in Proposition 5.9 $\pi$ by $R$ we are done.
	\hfill$\blacksquare$

	\vs{0.5cm}
	{\bf 5.11 Comment}
	
	\begin{itemize}
		\item[(1)]
		It is conjectured that $1/\pi$ is not a period. It is an integrated algebraic number.
		\item[(2)] 
		It is conjectured that the Euler number $e$ is not a period. We have that $e$ is an integrated algebraic number (see below).
		\item[(3)] 
		It is conjectured that the Euler constant $\gamma$ is not a period. We conjecture that it is not an integrated algebraic number. It is an exponential integrated algebraic number (see below).
	\end{itemize}
	
	\vs{0.2cm}
	{\bf 5.12 Proposition}
	{\it \begin{itemize}
			\item[(1)] We have $e\in \II\IA$.
			\item[(2)] We have $\gamma\in \II\IA_\exp$.
	\end{itemize}}
	{\bf Proof:}
	
	\vs{0.1cm}
	(1): By Theorem 2.2 and Proposition 3.13 we have $e=\exp(1)\in \II\IA$.
	
	\vs{0.2cm}
	(2): Of the many integral representations of $\gamma$ we use the following. We have (see for example Havil [13, p. 109])
	$$\gamma=\int_0^\infty\exp(-x)\Big(\frac{1}{1-\exp(-x)}-\frac{1}{x}\Big)\,dx.$$
	The integrand is obviously definable in $\widetilde{\ma{IA}}(\exp)$.
	Hence $\gamma\in \mathrm{Vol}(\widetilde{\ma{IA}}(\exp))=\II\IA_\exp$.
	\hfill$\blacksquare$
	
	\vs{0.5cm}
	By above we have the following picture:
	
	\vs{0.5cm}
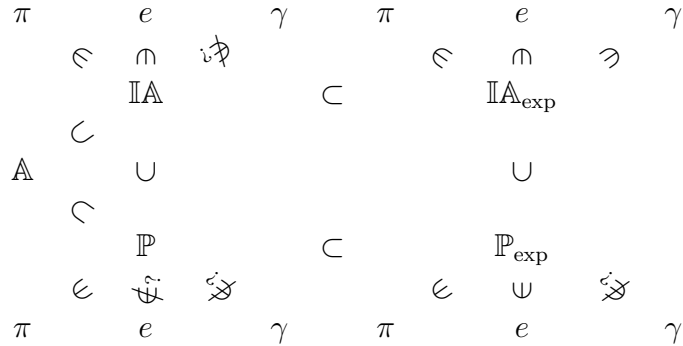
	\begin{center}
		\begin{tabular}{ccccccccccc}
			$\pi$&&$e$&&$\gamma$&&$\pi$&&$e$&&$\gamma$\\
			&\rotatebox[origin=c]{-35}{$\in$}&\rotatebox[origin=c]{-90}{$\in$}&\rotatebox[origin=c]{-145}{$\notin_?$}&&&&
			\rotatebox[origin=c]{-35}{$\in$}&\rotatebox[origin=c]{-90}{$\in$}&
			\rotatebox[origin=c]{-145}{$\in$}&\\
			&&$\II\IA$&&&$\subset$&&&$\II\IA_\exp$&&\\
			&\rotatebox[origin=c]{35}{$\subset$}&&&&&&&&&\\
			$\IA$&&\rotatebox[origin=c]{90}{$\subset$}&&&&&&\rotatebox[origin=c]{90}{$\subset$}&&\\
			&\rotatebox[origin=c]{-35}{$\subset$}&&&&&&&&&\\
			&&$\mathbb{P}$&&&$\subset$&&&$\mathbb{P}_\exp$&&\\
			&\rotatebox[origin=c]{35}{$\in$}&\rotatebox[origin=c]{90}{$\notin_?$}&\rotatebox[origin=c]{145}{$\notin_?$}&&&&\rotatebox[origin=c]{35}{$\in$}&
			\rotatebox[origin=c]{90}{$\in$}&
			\rotatebox[origin=c]{145}{$\notin_?$}&\\
			$\pi$&&$e$&&$\gamma$&&$\pi$&&$e$&&$\gamma$\\
		\end{tabular}
	\end{center}
	
	\vs{0.5cm}
	{\bf 5.13 Comment}
	
	\vs{0.1cm}
	By the above results the exponential integrated algebraic numbers are a good candidate for a natural class of real numbers beyond the periods capturing the  important mathematical constants as formulated by Kontsevich and Zagier [25]. Another candidate might also be the structure of the gobal exponential function and the restricted sine function on the reals (see van den Dries [8]) respectively its prime model. But we do not know whether there definable functions have good description in a suitable language and we do not have information about the universe of the prime model.

	\vs{2cm}
	{\bf References}
	
	{\footnotesize
		\begin{itemize}
			\item[(1)] 
			B. Bakker, Y. Brunebarbe, J. Tsimerman:
			o-minimal GAGA and a conjecture of Griffiths.  	
			{\it Invent. Math.} {\bf 232} (2023), no. 1, 163-228.
			
			\item[(2)] J. Bochnak, M. Coste, M.-F. Roy: Real algebraic geometry. Ergebnisse der Mathematik und ihrer Grenzgebiete {\bf 36}, Springer, 1998.
			
			\item[(3)] R. Cluckers, D. Miller: Stability under integration of sums of products of real globally subanalytic functions and their logarithms.
			{\it Duke Math. J.} {\bf 156} (2011), no. 2, 311-348.
			
			\item[(4)] J. Commelin, P. Habegger, A. Huber: 
			Exponential periods and o-minimality. arXiv:2007.08280
			
			\item[(5)] G. Comte, J.-M. Lion, J.-P. Rolin: Nature
			log-analytique du volume des sous-analytiques. {\it Illinois J.
				Math.} {\bf 44} (2000), no. 4, 884-888.
			
			\item[(6)] J. Denef, L. van den Dries: $p$-adic and real subanalytic sets. {\it Ann. of Math. (2)} {\bf 128} (1988), no. 1, 79-138.
			
			\item[(7)] J. Denef, L. Lipshitz: Ultraproducts and approximation in local rings. II. {\it Math. Ann.} {\bf 253} (1980), no.1, 1-28.
			
			\item[(8)] L. van den Dries: On the elementary theory of restricted elementary functions.
			{\it J. Symbolic Logic} {\bf 53}(1988), no.3 , 796-808.
			
			\item[(9)] L. van den Dries: Tame Topology and O-minimal structures. London Math. Soc. Lecture Note Ser. {\bf 248}, Cambridge University Press, 1998.
			
			\item[(10)] L. van den Dries. D. Marker, A. Macintyre: The elementary theory of restricted analytic fields with exponentiation. {\it Ann. of Math. (2)} {\bf 140} (1994), no. 1, 183-205.
			
			\item[(11)] L. van den Dries, P. Speissegger: O-minimal preparation theorems, pp. 87-116 in Model Theory and Applications,
			Quad. Mat. {\bf 11}, Aracne, Rome, 2002.
			
			\item[(12)] R. Gunning, H. Rossi: Analytic functions of several complex variables.
			Reprint of the 1965 original. AMS Chelsea Publishing, 2009.
			
			\item[(13)] L. Havil: Gamma - Exploring Euler's Constant.
			Princeton Sci. Lib.
			Princeton University Press, 2009.
			
			\item[(14)] R. Hardt, P. Lambrechts, V. Turchin, I. Voli\'c:
			Real homotopy theory of semi-algebraic sets. 
			{\it Algebr. Geom. Topol.} {\bf 11} (2011), no. 5, 2477-2545. 
			
			\item[(15)] A. Huber, S. M\"uller-Stach: Periods and Nori motives. With contributions by Benjamin Friedrich and Jonas von Wangenheim. Ergebnisse der Mathematik und ihrer Grenzgebiete. 3. Folge. A Series of Modern Surveys in Mathematics {\bf 65}. Springer, 2017.
			
			\item[(16)] A. Huber, G. W\"ustholz: Transcendence and linear relations of $1$-periods. Cambridge Tracts in Mathematics {\bf 227}. Cambridge University Press, 2022.
			
			\item[(17)] T. Kaiser: On convergence of integrals in
			o-minimal structures on archimedean real closed fields.
			{\it Annales Polonici Mathematici} {\bf 87} (2005), 175-192.
			
			\item[(18)]
			T. Kaiser: First order tameness of measures. 
			{\it Ann. Pure Appl. Logic} {\bf 163} (2012), no. 12, 1903-1927.
			
			\item[(19)] T. Kaiser: Integration of semialgebraic functions and integrated Nash functions.
			{\it Math. Zeitschrift} {\bf 275} (2013), no. 1-2, 349-366.
			
			\item[(20)] T. Kaiser: Multivariate Puiseux rings induced by a Weierstra\ss $ $ system and twisted group rings.
			{\it Communications in Algebra} {\bf 42} (2014), no. 11, 4619-4634.
			
			\item[(21)] T. Kaiser: Lebesgue measure and integration theory on non-archimedean real closed fields with archimedean value group.
			{\it Proceedings of the London Mathematical Society} {\bf 116} (2018), no. 2, 209-247.
			
			\item[(22)] T. Kaiser: Logarithms, constructible functions and integration on non-archimedean models of the theory of the real field with restricted analytic functions with value group of finite archimedean rank.
			{\it Fundamenta Mathematicae} {\bf 256} (2022), 285-306.
			
			\item[(23)] I. Kaplansky:
			Commutative Rings. Allen and Bacon,  1970.
			
			\item[(24)] M. Kontsevich, Y. Soibelman: Deformations of algebras over operads and the Deligne conjecture. Conf\'erence Mosh\'e Flato 1999, Vol. I (Dijon), 255-307, Math. Phys. Stud., {\bf 21}, Kluwer Acad. Publ., Dordrecht, 2000.
			
			\item[(25)]
			M. Kontsevich, D. Zagier:
			Periods. Mathematics unlimited -- 2001 and beyond, 771-808, Springer, 2001.
			
			\item[(26)] J.-M. Lion, J.-P. Rolin: Th\'{e}or\`{e}me de
			pr\'{e}paration pour les fonctions logarithmico-exponentielles.
			{\it Ann. Inst. Fourier} {\bf 47} (1997), no.3, 859-884.
			
			\item[(27)] J.-M. Lion, J.-P. Rolin: Int\'egration des fonctions
			sous-analytiques et volumes des sous-ensembles sous-analytiques.
			{\it Ann. Inst. Fourier} {\bf 48} (1998), no.3, 755-767.
			
			\item[(28)] D. Miller: A preparation theorem for Weierstrass systems.
			{\it Trans. Amer. Math. Soc.} {\bf 358} (2006), no.10, 4395-4439.
			
			\item[(29)] S. M\"uller-Stach: What is \ldots a period? {\it Notices Amer. Math. Soc.} {\bf 61} (2014), no. 8, 898-899.
			
			\item[(30)] A. Rambaud: Quasi-analycit\'e, o-minimalit\'e et \'elimination des quantificateurs. Ph.D. thesis, Universit\'e Paris 7 - Denis Diderot, 2005.
			
			\item[(31)] J. Ruiz: The Basic Theory of Power Series.
			Advanced Lectures in Mathematics. Vieweg, 1993.
			
			\item[(32)] K. Tent, M. Ziegler: Computable functions on the reals. {\it M\"unster J. of Math.} {\bf 3} (2010), 43-66.
	\end{itemize}}

	\vspace{1.5cm}
	\noi
	Tobias Kaiser \\
	Faculty of Computer Science and Mathematics\\
	University of Passau\\
	94030 Passau\\
	Germany\\
	email: tobias.kaiser@uni-passau.de
\end{document}